%% file: bic.tex
\documentclass[reqno]{amsart}

\usepackage[mathscr]{eucal}
\usepackage{amsmath}
\usepackage{amssymb}
\usepackage{latexsym}
\usepackage{amsthm}

\numberwithin{equation}{section}

\theoremstyle{plain}
\newtheorem{thm}{Theorem}[section]
\newtheorem{cor}[thm]{Corollary}
\newtheorem{prop}[thm]{Proposition}
\newtheorem{lem}[thm]{Lemma}
\newtheorem{mthm}{Main Theorem} %

\theoremstyle{definition}
\newtheorem{defn}[thm]{Definition}

\def\ti#1{\mbox{\tiny $#1$}}
\def\scri#1{\mbox{\scriptsize $#1$}}
\def\ft#1{\mbox{\footnotesize $#1$}}

\def\bol#1{\boldsymbol{#1}}

\newcommand{\al}{\alpha}
\newcommand{\be}{\beta}
\newcommand{\ga}{\gamma}
\newcommand{\de}{\delta}
\newcommand{\vep}{\varepsilon}
\newcommand{\lam}{\lambda}
\newcommand{\Lam}{\Lambda}
\newcommand{\ol}{\overline}
\newcommand{\wt}{\widetilde}

\newcommand{\til}{\tilde}
\newcommand{\st}{\stackrel}
\newcommand{\cl}{\mathcal} %
\newcommand{\fr}{\mathfrak} %
\newcommand{\bd}{\mathbf} %

\newcommand{\SI}{\st{\circ}{\bd I}}
\newcommand{\tiJ}{\ti{\bd J}}
\newcommand{\tiK}{\ti{\bd K}}
\newcommand{\tiL}{\ti{\bd L}}

\newcommand{\scK}{\scri{\bd K}}

\newcommand{\sPi}{\st{\circ}{\Pi}}
\newcommand{\sDe}{\st{\circ}{\De}}
\newcommand{\sW}{\st{\circ}{W}}
\newcommand{\sQ}{\st{\circ}{Q}}
\newcommand{\De}{\Delta}

\newcommand{\preq}{\preceq}
\newcommand{\re}{{\De}^{re}}
\newcommand{\pre}{{\De}^{re}_+}
\newcommand{\nre}{{\De}^{re}_-}
\newcommand{\im}{{\De}^{im}}
\newcommand{\pim}{{\De}^{im}_+}

\newcommand{\Lra}{\Longrightarrow}
\newcommand{\lra}{\longrightarrow}
\newcommand{\la}{\langle}
\newcommand{\ra}{\rangle}
\newcommand{\Ome}{\Omega}

\newcommand{\vPhi}{\varPhi}
\newcommand{\nab}{\nabla}

\newcommand{\bs}{\mbox{\boldmath $s$}}

\begin{document}
\title[Parametrizations of infinite biconvex sets]
{Parametrizations of infinite biconvex sets \\
in affine root systems}
\maketitle

\begin{center}
Ken Ito \\
\vspace{0.5em}
Graduate School of Mathematics, Nagoya University, \\
Chikusa-ku, Nagoya 464--8602, Japan \\
E-mail: ken-ito{\char'100}math.nagoya-u.ac.jp
\end{center}

\vspace{0.5em}
\begin{quotation}
\noindent
{\bf Abstract.}
 We investigate in detail relationships between the set
 ${\mathfrak B}^\infty$ of all infinite ``biconvex'' sets
 in the positive root system $\Delta_+$ of an arbitrary
 untwisted affine Lie algebra ${\mathfrak g}$ and
 the set ${\mathcal W}^\infty$ of all infinite ``reduced word''
 of the Weyl group of ${\mathfrak g}$.
 The study is applied to the classification of ``convex orders'' on
 $\Delta_+$ (cf. \cite{kI}), which are indispensable to construct
 ``convex bases'' of Poincar\'e-Birkhoff-Witt type of the upper
 triangular subalgebra $U_q^+$ of the quantized universal enveloping
 algebra $U_q({\mathfrak g})$.
 We construct a set $\boldsymbol{\mathcal P}$ by using data of
 the underlying finite-dimensional simple Lie algebra,
 and bijective mappings
 $\nabla\colon\boldsymbol{\mathcal P}\to{\mathfrak B}^\infty$ and
 $\chi\colon\boldsymbol{\mathcal P}\to W^\infty$ such that
 $\nabla=\varPhi^\infty\circ\chi$, where
 $W^\infty$ is an quotient set of ${\mathcal W}^\infty$ and
 $\varPhi^\infty\colon W^\infty\to{\mathfrak B}^\infty$
 is a natural injective mapping.

\noindent
{\bf Mathematics Subject Classifications (2000).}
 20F55, 17B37.

\noindent
{\bf Key Words.}
 convex order, convex basis, biconvex set, reduced word,
 affine root system, Weyl group, Coxeter group,
 Kac-Moody Lie algebra, quantum algebra.
\end{quotation}

\vspace{0.5em}
\input{sec1}
\input{sec2}
\input{sec3}
\input{sec4}
\input{sec5}
\input{sec6}
\input{sec7}

\vspace{0.5em}
\begin{center}
{\sc Acknowledgements}
\end{center}

 The author would like to thank Prof. Akihiro Tsuchiya and
 Prof. Takahiro Hayashi for constant help and precious advice.

\end{document}

%% file: sec1.tex
\section{Introduction}
 Let $\De$ be the root system of a Kac-Moody Lie algebra
 ${\fr g}$, $\De_+$ (resp. $\De_-$) the set of all positive
 (resp. negative) roots relative to the root basis
 $\Pi=\{\al_i\,|\,i\in{\bd I}\}$, and
 $W=\la{\,s_i\,|\,i\in{\bd I}\,}\ra$ the Weyl group of ${\fr g}$,
 where $s_i$ is the reflection associated with $\al_i$.
 Then $(W,S)$ is a Coxeter system with
 $S=\{\,s_i\,|\,i\in{\bd I}\,\}$ (cf. \cite{vK}).
 We call an infinite sequence
 $\bs=(\bs\ft{(p)})_{p\in{\mathbb N}}\in S^{\mathbb N}$
 an {\it infinite reduced word\/} of $(W,S)$
 if the length of the element
 $z_{\bs}\ft{(p)}:=\bs\ft{(1)}\cdots\bs\ft{(p)}\in W$ is $p$
 for each $p\in{\mathbb N}$, and call a subset $B\subset\De_+$
 a {\it biconvex set\/} if it satisfies the following conditions:

 C(i)\quad %
 $\be,\,\ga\in B,\;\; \be+\ga\in\De_+\,\Lra\,\be+\ga\in B$;

 C(ii)\quad %
 $\be,\,\ga\in\De_+\setminus B,\;\; \be+\ga\in\De_+ %
 \,\Lra\, \be+\ga\in\De_+\setminus B$.

\noindent
 If, in addition, $B$ is a subset of the set $\pre$ of
 all positive real roots, then $B$ is called
 a {\it real biconvex set\/}. The purpose of this paper is
 to investigate in detail relationships between infinite reduced
 words and infinite real biconvex sets in the case where ${\fr g}$
 is an arbitrary untwisted affine Lie algebra.

 Before explaining the detail of our work, we will explain
 the background of the theory of infinite reduced words and
 infinite real biconvex sets.
 The motive of this study is related to the construction of
 {\it convex bases\/} of the upper triangular subalgebra $U_q^+$
 of the quantized universal enveloping algebra $U_q({\fr g})$.
 Convex bases are Poincar\'e-Birkhoff-Witt type bases
 with a convex property concerning the ``$q$-commutator''
 of two ``$q$-root vectors'' of $U_q^+$.
 The convex property is useful for calculating values of
 the standard bilinear pairing between $U_q^+$ and
 the lower triangular subalgebra $U_q^-$, and is applied to
 explicit calculations of the universal $R$-matrix of
 $U_q({\fr g})$ (cf. \cite{KR}, \cite{KT}).
 By the way, each convex basis of $U_q^+$ is formed by monomials
 in certain $q$-root vectors of $U_q^+$ multiplied in
 a predetermined total order with a convex property on $\De_+$.
 Such a total order on $\De_+$ called a {\it convex order\/}.

 In the case where ${\fr g}$ is an arbitrary finite-dimensional
 simple Lie algebra, it is known that there exists a natural
 bijective mapping from the set of all convex orders on $\De_+$ to
 the set of all reduced expressions of the longest element of $W$
 (cf. \cite{pP}), and G.~Lusztig constructed convex bases of $U_q^+$
 associated with all reduced expressions of the longest element of $W$
 by using a braid group action on $U_q$ (cf. \cite{gL}).
 Thus all convex bases of $U_q^+$ can be constructed in this case.

 In the case where ${\fr g}$ is an arbitrary untwisted affine
 Lie algebra, convex orders on $\De_+$ are closely related to
 infinite reduced words of $(W,S)$.
 More precisely, each infinite reduced word naturally
 corresponds to a ``{\it 1-row type\/}'' convex order on
 an infinite real biconvex set.
 Moreover, each convex order on $\De_+$ is made from each couple of
 ``{\it maximal\/}'' (infinite) real biconvex sets with
 convex orders which divides $\pre$ into two parts (cf. \cite{kI}).
 In \cite{jB}, J. Beck constructed convex bases of $U_q^+$
 associated with convex orders on $\De_+$ arisen from a certain
 couple of maximal real biconvex sets with 1-row type convex orders
 which divides $\pre$ into two parts.
 However, we expect that it is possible to generalize Beck's
 construction, since we find that there exist several types of
 convex orders called ``{\it n-row types\/}'' on each maximal
 real biconvex set which are not used in Beck's  construction
 (cf. \cite{kI}).
 To analyze convex orders on maximal real biconvex sets,
 it is important to consider the following two problems:
 (1) classify all infinite real biconvex sets;
 (2) describe in detail relationships between
 the set of all infinite reduced words and
 the set of all infinite real biconvex sets.
 In this paper, we concentrate on the above two problems
 for the untwisted affine cases. By using results of this paper,
 in \cite{kI}, we classified all convex orders on $\De_+$, and then
 gave a general method of constructing convex orders on $\De_+$
 for the untwisted affine cases.

 This paper is organized as follows.
 In Section 2, we first give the definition of
 biconvex sets for a class of root systems with
 Coxeter group actions, and then state several fundamental results.
 We next define infinite reduced words for each
 Coxeter system $(W,S)$ and an equivalence relation $\sim$
 on the set ${\cl W}^\infty$ of all infinite reduced words,
 and then define $W^\infty$ to be the quotient set
 of ${\cl W}^\infty$ relative to $\sim$.
 We next define an injective mapping
 $\vPhi^\infty\colon W^\infty\to{\fr B}^\infty$,
 where ${\fr B}^\infty$ is the set of
 all infinite real biconvex sets.
 At the end of Section 2, we define a left action
 of $W$ on $W^\infty$ which plays an important role
 in the proof of the main theorem.
 In Section 3,
 we introduce notation for the untwisted affine cases.
 In Section 4,
 we give preliminary results for classical root systems.
 From Section 5 to Section 7,
 we treat only the untwisted affine cases.
 In Section 5,
 we give several methods of constructing biconvex sets.
 In Section 6,
 we give a Parametrization of real biconvex sets.
 In Section 7,
 we give the following main results.

\begin{mthm}
 If ${\fr g}$ is an arbitrary untwisted affine Lie algebra,
 there exist Parametrizations {\em(bijective mappings)}
 $\nab\colon\bol{\cl P}\to{\fr B}^\infty$ and
 $\chi\colon\bol{\cl P}\to W^\infty$ such that
 the following diagram is commutative\/{\em:}

\setlength{\unitlength}{0.5mm}
\begin{picture}(240,50)(-100,-5)
\put(10,30){\makebox(20,10){${\fr B}^\infty$}}
\put(-20,0){\makebox(20,10){$W^\infty$}}
\put(42,0){\makebox(20,10){$\bol{\cl P}$\;,}}
\put(45,5){\vector(-1,0){48}}
\put(-7,10){\vector(1,1){20}}
\put(47,10){\vector(-1,1){20}}
\put(17,0){$\chi$}
\put(-7,20){$\vPhi^\infty$}
\put(40,20){$\nab$}
%
\end{picture}

\noindent
 where the set $\bol{\cl P}$ is defined by using data of
 the underlying finite-dimensional simple Lie algebra
 $\st{\circ}{\fr g}$ {\em(see Definition 6.5)}.
 In particular, $\vPhi^\infty$ is bijective. Moreover,
 $W^\infty$ decomposes into the direct finite sum of
 orbits relative to a left action of $W$.
\end{mthm}

 Note that in \cite{CP} P.~Cellini and P.~Papi showed that
 if $B$ is an infinite real biconvex set,
 then there exist $v,\,t\in W$ such that $t$ is a translation,
 $\ell(vt)=\ell(v)+\ell(t)$, and $B=\cup_{k\geq0}\vPhi(vt^k)$,
 where $\vPhi(z)=\{\,\be\in\De_+\;|\;z^{-1}(\be)\in\De_-\,\}$.

%% file: sec2.tex
\section{Definitions and several results}
 Let ${\mathbb R}$, ${\mathbb Z}$, and ${\mathbb N}$
 be the set of the real numbers, the integers,
 and the positive integers, respectively.
 Set
 ${\mathbb R}_{\geq a}:=\{b\in{\mathbb R}\,|\,b\geq a\}$,
 ${\mathbb R}_{>a}:=\{b\in{\mathbb R}\,|\,b>a\}$,
 and ${\mathbb Z}_{\geq a}:={\mathbb Z}\cap{\mathbb R}_{\geq a}$
 for each $a\in{\mathbb R}$.
 Denote by $\sharp{U}$ the cardinality of a set $U$.
 When $A$ and $B$ are subsets of $U$, write
 $A\,\dot{\subset}\,B$ or $B\,\dot{\supset}\,A$
 if $\sharp(A\setminus B)<\infty$, and write $A\,\dot{=}\,B$
 if both $A\,\dot{\subset}\,B$ and $A\,\dot{\supset}\,B$.
 Then $\dot{=}$ is an equivalence relation on
 the power set of $U$.

\vspace{0.5em}
 Let $W$ be a group generated by a set $S$ of
 involutive generators (i.e., $s\neq1,\,s^2=1,\;\forall\,s\in S$),
 and $(V,\De,\Pi)$ a triplet satisfying
 the following conditions R(i)--R(iv).

 R(i)\; It consists of a representation space $V$ of $W$ over
 ${\mathbb R}$, a $W$-invariant subset $\De\subset V\setminus\{0\}$
 which is symmetric (i.e., $\De=-\De$), and
 a subset $\Pi=\{\,\al_s\,|\,s\in S\,\}\subset\De$.

 R(ii)\; Each element of $\De$ can be written as
 $\sum_{s\in S}a_s\al_s$ with either $a_s\geq0$ for all $s\in S$
 or $a_s\leq0$ for all $s\in S$, but not in both ways.
 Accordingly, we write $\al>0$ or $\al<0$, and put
 $\De_+=\{\,\al\in\De\,|\,\al>0\,\}$ and
 $\De_-=\{\,\al\in\De\,|\,\al<0\,\}$.

 R(iii)\; For each $s\in S$, $s(\al_s)=-\al_s$ and
 $s(\De_+\setminus\{\al_s\})=\De_+\setminus\{\al_s\}$.

 R(iv)\; If $w\in W$ and $s,\,s'\in S$ satisfy
 $w(\al_{s'})=\al_s$, then $ws'w^{-1}=s$.

\begin{defn} 
 Define subsets $\re$, $\im$, $\re_\pm$,
 and $\im_\pm$ of $\De$ by setting
\begin{align*}
 \re&:=\{\,w(\al_s)\;|\;w\in W,\,s\in S\,\},\qquad %
 \im:=\De\setminus\re, \\
 \re_\pm&:=\re\cap\De_\pm,\qquad\im_\pm:=\im\cap\De_\pm.
\end{align*}
 Note that $W$ stabilizes $\re$ and $\pim$.
 We also set
\[ \vPhi(y):=\{\,\be\in\De_+\;|\;y^{-1}(\be)<0\,\} \]
 for each $y\in W$. Note that $\vPhi(y)\subset\pre$.
\end{defn}

\begin{thm}[\cite{vD2}] 
 The pair $(W,S)$ is a Coxeter system, i.e., it satisfies
 the exchange condition. Moreover, if $y=s_1s_2\cdots s_n$
 with $n\in{\mathbb N}$ and $s_1,s_2,\dots,s_n\in S$ is
 a reduced expression of an element $y\in W\setminus\{1\}$, then
\[ \vPhi(y)=\{\,\al_{s_1},\,s_1(\al_{s_2}),
 \,\dots,\,s_1\cdots s_{n-1}(\al_{s_n})\,\} \]
 and the elements of $\vPhi(y)$ displayed above are distinct
 from each other. In particular, $\sharp{\vPhi(y)}=\ell(y)$,
 where $\ell\colon W\to{\mathbb Z}_{\geq0}$ is the length function
 of $(W,S)$.
\end{thm}

\noindent
{\it Remarks.}
 (1) For each Coxeter system $(W,S)$, a triplet $(V,\De,\Pi)$
 is called a {\it root system\/} of $(W,S)$ if
 it satisfies the conditions R(i)--R(iv).

 (2) Let $\sigma\colon W\to\mathrm{GL}(V)$ be
 the geometric representation of a Coxeter system $(W,S)$
 (cf. \cite{nB}), where $V$ is a real vector space with
 a basis $\Pi=\{\,\al_s\,|\,s\in S\,\}$.
 Then $(V,\De,\Pi)$ is a root system of $(W,S)$
 (cf. \cite{vD2}), where
 $\De=\{\,\sigma(w)(\al_s)\,|\,w\in W,\,s\in S\,\}$.
 We call it the root system associated with
 the geometric representation.

 (3) Let ${\fr g}$ be a Kac-Moody Lie algebra over ${\mathbb R}$
 with a Cartan subalgebra ${\fr h}$,
 $\De\subset{\fr h}^*\setminus\{0\}$ the root system of ${\fr g}$,
 $\Pi=\{\,\al_i\,|\,i\in{\bd I}\,\}$ a root basis of $\De$, and
 ${W}=\la{\,s_i\,|\,i\in{\bd I}\,}\ra\subset\mathrm{GL}({\fr h}^*)$
 the Weyl group of ${\fr g}$, where ${\fr h}^*$ is the dual space of
 ${\fr h}$ and $s_i$ is the simple reflection associated with $\al_i$
 (cf. \cite{vK}).
 Then $({\fr h}^*,\De,\Pi)$ is a root system of
 a Coxeter system $(W,S)$, where $S=\{\,s_i\,|\,i\in{\bd I}\,\}$.

\begin{lem} 
 Let $y_1$ and $y_2$ be elements of $W$.
\begin{itemize}
\item[(1)] We have
 $\vPhi(y_1y_2)\setminus{\vPhi(y_1)}\subset{y_1\vPhi(y_2)}$.
\item[(2)] If $y_1\vPhi(y_2)\subset\De_+$,
 then $\vPhi(y_1)\amalg y_1\vPhi(y_2)=\vPhi(y_1y_2)$.
\item[(3)] If $\vPhi(y_1)\subset{\vPhi(y_2)}$,
 then $\vPhi(y_2)=\vPhi(y_1)\amalg{y_1\vPhi(y_1^{-1}y_2)}$.
\item[(4)] The following two conditions are equivalent\/{\em:}
\[ \text{\em(i)}\;\; \ell(y_2)-\ell(y_1)=\ell(y_1^{-1}y_2);\qquad
 \text{\em(ii)}\;\; \vPhi(y_1)\subset\vPhi(y_2). \]
\end{itemize}
\end{lem}

\begin{proof}
 (1) Suppose that $\be\in{\vPhi(y_1y_2)}\setminus{\vPhi(y_1)}$.
 Then we have $y_1^{-1}(\be)>0$ and $y_2^{-1}(y_1^{-1}(\be))<0$.
 Thus we get
 $y_1^{-1}(\be)\in{\vPhi(y_2)}$ or $\be\in{y_1\vPhi(y_2)}$.

 (2) If $\be\in{y_1\vPhi(y_2)}$ then $y_1^{-1}(\be)>0$,
 and hence $\be\notin{\vPhi(y_1)}$.
 Thus we get $\vPhi(y_1)\cap{y_1\vPhi(y_2)}=\emptyset$.
 Hence, by (1) we have
 $\vPhi(y_1y_2)\subset{\vPhi(y_1)}\amalg{y_1\vPhi(y_2)}$.
 We next prove that $\vPhi(y_1)\subset{\vPhi(y_1y_2)}$.
 Suppose that $\be\in{\vPhi(y_1)}$ satisfies
 $\be\notin{\vPhi(y_1y_2)}$.
 Then we have $y_1^{-1}(\be)<0$ and $y_2^{-1}(y_1^{-1}(\be))>0$,
 which imply $-y_1^{-1}(\be)\in{\vPhi(y_2)}$.
 This contradicts to the assumption.
 Thus we get $\vPhi(y_1)\subset{\vPhi(y_1y_2)}$.
 We next prove that $y_1\vPhi(y_2)\subset{\vPhi(y_1y_2)}$.
 If $\be\in{y_1\vPhi(y_2)}$ then $y_1^{-1}(\be)\in{\vPhi(y_2)}$,
 and hence $y_2^{-1}(y_1^{-1}(\be))<0$.
 Thus we get $\be\in{\vPhi(y_1y_2)}$.
 Therefore
 $\vPhi(y_1)\amalg{y_1\vPhi(y_2)}\subset{\vPhi(y_1y_2)}$.

 (3) We first prove that $y_1\vPhi(y_1^{-1}y_2)\subset\De_+$.
 Suppose that $\be\in\vPhi(y_1^{-1}y_2)$ satisfies $y_1(\be)<0$.
 Then we have $-y_1(\be)\in\vPhi(y_1)\subset\vPhi(y_2)$,
 which implies $y_2^{-1}y_1(\be)>0$.
 This contradicts to $\be\in \vPhi(y_1^{-1}y_2)$.
 Thus we get $y_1\vPhi(y_1^{-1}y_2)\subset\De_+$, and hence
 $\vPhi(y_2)=\vPhi(y_1)\amalg y_1\vPhi(y_1^{-1}y_2)$ by (2).

 (4)(i)$\Rightarrow$(ii)
 By Theorem 2.2, we have
\begin{align*}
 \ell(y_2)-\ell(y_1) &\leq %
 \sharp y_1^{-1}\{\vPhi(y_2)\setminus{\vPhi(y_1)}\} \\
 &\leq \sharp \vPhi(y_1^{-1}y_2)= %
 \ell(y_1^{-1}y_2)=\ell(y_2)-\ell(y_1),
\end{align*}
 where the second inequality follows from (1).
 Thus we get
 $\sharp y_1^{-1}\{\vPhi(y_2)\setminus\vPhi(y_1)\}= %
 \ell(y_2)-\ell(y_1)$, and hence
 $\vPhi(y_1)\subset \vPhi(y_2)$.

 (ii)$\Rightarrow$(i)
 By (3) and Theorem 2.2, we get
 $\ell(y_2)=\ell(y_1)+\ell(y_1^{-1}y_2)$.
\end{proof}

\begin{defn} 
 For subsets $A,\,B\subset\De_+$ satisfying $B\subset A$,
 we call $B$ a {\it convex set in $A$\/}
 if it satisfies the following condition:

 $\mbox{C(i)}_{\ti{A}}\quad
 \be,\,\ga\in B,\;\; \be+\ga\in A \,\Lra\, \be+\ga\in B$.

\noindent
 We also call $B$ a {\it coconvex set in $A$\/}
 if it satisfies the following condition:

 $\mbox{C(ii)}_{\ti{A}}\quad
 \be,\,\ga\in A\setminus B,\;\; \be+\ga\in A %
 \,\Lra\, \be+\ga\in A\setminus B$.

\noindent
 Note that $B$ is a coconvex set in $A$ if and only if
 $A\setminus B$ is a convex set in $A$. Furthermore,
 we call $B$ a {\it biconvex set in $A$\/} if $B$
 is both a convex set in $A$ and a coconvex set in $A$.
 If, in addition, $B\subset\pre$, then $B$ is said to be
 a {\it real convex set in $A$}, a {\it real coconvex set in $A$},
 or a {\it real biconvex set in $A$} if $B$ is a convex set in $A$,
 a coconvex set in $A$, or a biconvex set in $A$, respectively.

 We will say simply that $B$ is
 a {\it convex set\/}, a {\it real convex set\/},
 a {\it coconvex set\/}, a {\it real coconvex set\/},
 a {\it biconvex set\/}, or a {\it real biconvex set\/} if $B$ is
 a convex set in $\De_+$, a real convex set in $\De_+$,
 a coconvex set in $\De_+$, a real coconvex set in $\De_+$,
 a biconvex set in $\De_+$, or a real biconvex set in $\De_+$,
 respectively.
 We denote $\mbox{C(i)}_{\ti{\De_+}}$ and $\mbox{C(ii)}_{\ti{\De_+}}$
 simply by C(i) and C(ii), respectively.
 We denote by ${\fr B}$ the set of all finite biconvex sets, and
 by ${\fr B}^\infty$ the set of all infinite real biconvex sets.
\end{defn}

\noindent
{\it Remark.}
 The condition $\mbox{C(ii)}_{\ti{A}}$ is equivalent to
 the following condition:
\[ \be,\,\ga\in A,\;\; \be+\ga\in B
 \,\Lra\, \be\in B \;\mbox{ or }\; \ga\in B. \]

 For each subsets $A,\,B\subset\De$, we set
\[ A\dotplus B:=\{\,\al+\be\;|\;\al\in A,\,\be\in B\,\}\cap\De. \]

\begin{lem} 
 Let $A$, $B$, and $C$ be subsets of $\De_+$ satisfying
 $B,\,C\subset A$, and $\{B_\lam\}_{\lam\in\Lam}$
 a family of subsets of $A$.
\begin{itemize}
\item[(1)] If $B$ is a biconvex set in $A$, then
 $A\setminus B$ is biconvex in $A$.
\item[(2)] If $B$ is a biconvex set in $A$, then
 $B\cap C$ is a biconvex set in $C$.
\item[(3)] Suppose that $B\subset C$ and
 $C$ is a convex set in $A$.
 Then $B$ is a convex set in $C$ if and only if
 $B$ is a convex set in $A$.
\item[(4)] If $(B_\lam\dotplus B_{\lam'})\cap A\subset
 \cup_{\lam\in\Lam}B_\lam$ for each $\lam,\,\lam'\in\Lam$,
 then $\cup_{\lam\in\Lam}B_\lam$ is a convex set in $A$.
\item[(5)] If $B_\lam$ is a convex set in $A$ for each $\lam\in\Lam$,
 then $\cap_{\lam\in\Lambda}B_{\lam}$ is a convex set in $A$.
\item[(6)] If $B_\lam$ is a biconvex set in $A$ for each $\lam$
 and $\preq$ is a total order on $\Lam$ such that
 $B_\lam\subseteq{B_{\lam'}}$ for each $\lam\preq\lam'$,
 then both $\cup_{\lam\in\Lam}B_\lam$ and
 $\cap_{\lam\in\Lam}B_\lam$ are biconvex sets in $A$.
\end{itemize}
\end{lem}

\begin{proof}
 (1)--(5) They are obvious.

 (6) Set $B_1:=\cup_{\lam\in\Lam}B_\lam$.
 By the assumption on the total order $\preq$,
 the family $\{B_\lam\}_{\lam\in\Lam}$ satisfies
 the sufficient condition in (4).
 Hence, $B_1$ is a convex set in $A$.
 On the other hand, since
 $A\setminus{B_1}=\cap_{\lam\in\Lam}(A\setminus{B_\lam})$,
 $A\setminus{B_1}$ is a convex set in $A$ by (1) and (5).
 Thus $B_1$ is a biconvex set in $A$.
 Set $B_2:=\cap_{\lam\in\Lam}B_\lam$.
 Let $\preq^{op}$ be the opposite order of $\preq$.
 Then $A\setminus{B_\lam}\subseteq{A}\setminus{B_{\lam'}}$
 if $\lam\preq^{op}\lam'$.
 Hence, $B_3:=\cup_{\lam\in\Lam}(A\setminus{B_\lam})$
 is a biconvex set in $A$.
 Thus $B_2$ is a biconvex set in $A$ since $B_2=A\setminus{B_3}$.
\end{proof}

\begin{thm}[\cite{pP}] 
 The assignment $y\mapsto\vPhi(y)$ defines an injective mapping
 from $W$ to ${\fr B}$. Moreover, if the root system
 $(V,\De,\Pi)$ satisfies the following two conditions then
 $\vPhi$ is surjective\/{\em:}

 {\em R(v)} each $\al\in\De_+\setminus\Pi$ can be written as
 $\be+\ga$ with $\be,\,\ga\in\De_+${\em;}

 {\em R(vi)} there exists a mapping
 $\mathrm{ht}\colon\De_+\to{\mathbb R}_{>0}$ such that
 $\mathrm{ht}(\be+\ga)=\mathrm{ht}(\be)+\mathrm{ht}(\ga)$ if
 $\be+\ga\in\De+$ for some $\be,\,\ga\in\De_+$.
\end{thm}

\noindent
{\it Remarks.}
 (1) The surjectivity of the mapping follows from the fact
 that if $C$ is a non-empty finite coconvex set then
 $C\cap\Pi\neq\emptyset$.
 The conditions R(v) and R(vi) are used to prove the fact.

 (2) Suppose that the root system $(V,\De,\Pi)$ satisfies
 the following two condition instead of R(v) and R(vi):

 $\mbox{R(v)}'$\quad
 each $\al\in\De_+\setminus\Pi$ can be written as $b\be+c\ga$
 with $b,\,c\in{\mathbb R}_{\geq1}$ and $\be,\,\ga\in\De_+$;

 $\mbox{R(vi)}'$\quad
 there exists a mapping
 $\mathrm{ht}\colon\De_+\to{\mathbb R}_{>0}$ such that
 $\mathrm{ht}(b\be+c\ga)=b\mathrm{ht}(\be)+c\mathrm{ht}(\ga)$ if
 $b\be+c\ga\in\De+$ for some $b,\,c\in{\mathbb R}_{>0}$ and
 $\be,\,\ga\in\De_+$.

\noindent
 Then $\vPhi$ is still surjective if ${\fr B}$ is replaced by
 the set of all finite subsets $B\subset\De_+$ satisfying
 the following two conditions:

 $\mbox{C(i)}'$\quad
 $\be,\,\ga\in B,\;\; b,\,c\in{\mathbb R}_{>0},\;\;
 b\be+c\ga\in\De_+ \,\Lra\, b\be+c\ga\in B$;

 $\mbox{C(ii)}'$\quad
 $\be,\,\ga\in\De_+\setminus B,\;\; b,\,c\in{\mathbb R}_{>0},\;\;
 b\be+c\ga\in\De_+ \,\Lra\, b\be+c\ga\in\De_+\setminus B$.

 (3) Let $(V,\De,\Pi)$ be the root system of a Coxeter system
 $(W,S)$ associated with the geometric representation.
 Then $(V,\De,\Pi)$ satisfies $\mbox{R(v)}'$ and $\mbox{R(vi)}'$.
 The condition $\mbox{R(v)}'$ is easily checked by
 reforming the proof of Proposition 2.1 in \cite{vD1}.
 Since $\Pi$ is linearly independent, we can define
 a mapping $\mathrm{ht}\colon\De_+\to{\mathbb R}_{>0}$
 by setting $\mathrm{ht}(\al):=\sum_{s\in S}a_s$
 for each $\al\in\De_+$, where $a_s$'s are non-negative
 real numbers such that $\al=\sum_{s\in S}a_s\al_s$.
 Then the mapping $\mathrm{ht}$ satisfies
 the required property in $\mbox{R(v)}'$.

\begin{defn} 
 For each infinite sequence
 $\bs=(\bs\ft{(p)})_{p\in{\mathbb N}}\in S^{\mathbb N}$,
 we define two mappings $z_{\bs}\colon{\mathbb N}\to W$ and
 $\phi_{\bs}\colon{\mathbb N}\to\re$ by setting
\[ z_{\bs}\ft{(p)}:=\bs\ft{(1)}\cdots\bs\ft{(p)},\qquad
 \phi_{\bs}\ft{(p)}:=z_{\bs}\ft{(p-1)}(\al_{\bs(p)}) \]
 for each $p\in{\mathbb N}$, where $z_{\bs}\ft{(0)}:=1$, and
 define a mapping $\wt{\vPhi}^\infty$ from $S^{\mathbb N}$ to
 the power set of $\pre$ by setting
\[ \wt{\vPhi}^\infty(\bs):=
 \bigcup_{p\in{\mathbb N}}\vPhi(z_{\bs}\ft{(p)}) \]
 for each $\bs\in S^{\mathbb N}$.
 We call an element $\bs\in S^{\mathbb N}$
 an {\it infinite reduced word\/} of $(W,S)$ if
 $\ell(z_{\bs}\ft{(p)})=p$ for all $p\in{\mathbb N}$,
 and denote by ${\cl W}^\infty$ the subset of $S^{\mathbb N}$
 of all infinite reduced words of $(W,S)$.
\end{defn}

\begin{lem} 
 For a pair $(\bs,\bs')$ of elements of ${\cl W}^\infty$, we write
 $\bs\sim\bs'$ if for each $(p,q)\in{\mathbb N}^2$ there exists
 $(p_0,q_0)\in{\mathbb Z}_{\geq p}\times{\mathbb Z}_{\geq q}$
 such that
\[ \ell(z_{\bs}\ft{(p)}^{-1}z_{\bs'}\ft{(p_0)})=p_0-p,\quad
 \ell(z_{\bs'}\ft{(q)}^{-1}z_{\bs}\ft{(q_0)})=q_0-q. \]
 Then $\sim$ is an equivalence relation on ${\cl W}^\infty$.
\end{lem}

\begin{proof}
 The reflexive row and the symmetric row are obvious.
 To prove the transitive row,
 suppose that $\bs\sim\bs',\,\bs'\sim\bs''$
 for some $\bs,\,\bs',\,\bs''\in{\cl W}^\infty$.
 For each $p\in{\mathbb N}$, choose $p_0\geq p$ and $p_1\geq p_0$
 satisfying
 $\ell(z_{\bs}\ft{(p)}^{-1}z_{\bs'}\ft{(p_0)})=p_0-p$ and
 $\ell(z_{\bs'}\ft{(p_0)}^{-1}z_{\bs''}\ft{(p_1)})=p_1-p_0$.
 Then we have
\begin{align*}
 p_1-p=|\ell(z_{\bs}\ft{(p)}^{-1})-\ell(z_{\bs''}\ft{(p_1)})|
 &\leq \ell(z_{\bs}\ft{(p)}^{-1}z_{\bs''}\ft{(p_1)}) \\
 &\leq \ell(z_{\bs}\ft{(p)}^{-1}z_{\bs'}\ft{(p_0)})+
 \ell(z_{\bs'}\ft{(p_0)}^{-1}z_{\bs''}\ft{(p_1)}) \\
 &= (p_0-p)+(p_1-p_0)=p_1-p.
\end{align*}
 Thus we get
 $\ell(z_{\bs}\ft{(p)}^{-1}z_{\bs''}\ft{(p_1)})=p_1-p$.
 Similarly, we see that for each $q\in{\mathbb N}$
 there exists $q_1\in{\mathbb Z}_{\geq q}$ such that
 $\ell(z_{\bs''}\ft{(q)}^{-1}z_{\bs}\ft{(q_1)})=q_1-q$.
 Therefore we get $\bs\sim\bs''$.
\end{proof}

\begin{defn} 
 We denote by $W^\infty$ the quotient set of ${\cl W}^\infty$
 relative to the equivalence relation $\sim$, and by $[\bs]$
 the coset containing $\bs\in{\cl W}^\infty$.
\end{defn}

\begin{prop} 
 Let $\bs$ and $\bs'$ be elements of $S^{\mathbb N}$.

 {\em(1)} We have $\bs\in{\cl W}^\infty$ if and only if
 $\phi_{\bs}\ft{(p)}>0$ for all $p\in{\mathbb N}$.

 {\em(2)} If $\bs\in{\cl W}^\infty$, then
 $\wt{\vPhi}^\infty(\bs)=
 \{\,\phi_{\bs}\ft{(p)}\;|\;p\in{\mathbb N}\,\}$
 and the all elements $\phi_{\bs}\ft{(p)}$ of
 $\wt{\vPhi}^\infty(\bs)$ are distinct from each other.

 {\em(3)} If $\bs\in{\cl W}^\infty$, then
 $\wt{\vPhi}^\infty(\bs)\in{\fr B}^\infty$.

 {\em(4)} Suppose that $\bs,\,\bs'\in{\cl W}^\infty$.
 Then $\bs\sim\bs'$ if and only if
 $\wt{\vPhi}^\infty(\bs)=\wt{\vPhi}^\infty(\bs')$.
\end{prop}

\begin{proof}
 (1) We see that $\bs\in{\cl W}^\infty$ if and only if
 $\ell(z_{\bs}\ft{(p-1)}\bs\ft{(p)})>\ell(z_{\bs}\ft{(p-1)})$
 for all $p\in{\mathbb N}$.
 Hence the assertion follows from the fact that
 $z(\al_s)>0$ if and only if $\ell(zs)>\ell(z)$
 for $z\in W$ and $s\in S$.

 (2) This follows from Theorem 2.2.

 (3) We see that $\wt{\vPhi}^\infty(\bs)$ is an infinite set
 by (2). For each $p\leq q$, we have
 $\vPhi(z_{\bs}\ft{(p)})\subseteq{\vPhi(z_{\bs}\ft{(q)})}$.
 Thus we get $\wt{\vPhi}^\infty(\bs)\in{\fr B}^\infty$
 by Lemma 2.5(6) and Theorem 2.6.

 (4) By Lemma 2.3(4), we see that the condition $\bs\sim\bs'$
 is equivalent to the condition that
 for each $(p,q)\in{\mathbb N}^2$ there exists
 $(p_0,q_0)\in{\mathbb Z}_{\geq p}\times{\mathbb Z}_{\geq q}$
 such that
 $\vPhi(z_{\bs}\ft{(p)})\subset{\vPhi(z_{\bs'}\ft{(p_0)})}$ and
 $\vPhi(z_{\bs'}\ft{(q)})\subset{\vPhi(z_{\bs}\ft{(q_0)})}$.
 Thus $\bs\sim\bs'$ if and only if
 $\wt{\vPhi}^\infty(\bs)=\wt{\vPhi}^\infty(\bs')$.
\end{proof}

\begin{defn} 
 Thanks to Proposition 2.10(3)(4), we have an injective mapping
\[ \vPhi^\infty\colon W^\infty\,\lra\,{\fr B}^\infty,\quad
 [\bs]\,\longmapsto\,\vPhi^\infty([\bs]):=\wt{\vPhi}^\infty(\bs). \]
\end{defn}

 We define a left action of $W$ on $W^\infty$.

\begin{defn} 
 For each $x\in W$ and $\bs\in S^{\mathbb N}$, we set
\[ \wt{\vPhi}^\infty(x,\bs):=\{\,\be\in\pre\;|\;
 \exists p_0\in{\mathbb N};\,\forall p\geq p_0,
 \,(xz_{\bs}\ft{(p)})^{-1}(\be)<0\,\}. \]
\end{defn}

\begin{lem} 
 {\em(1)} If $\bs\in{\cl W}^\infty$, then
 $\wt{\vPhi}^\infty(1,\bs)=\wt{\vPhi}^\infty(\bs)$.

 {\em(2)} If $x\in W$ and $\bs\in{\cl W}^\infty$, then
 there exists an element $\bs'\in{\cl W}^\infty$ such that
 $\wt{\vPhi}^\infty(\bs')=\wt{\vPhi}^\infty(x,\bs)$.
 More precisely, a required $\bs'$ can be constructed by
 applying the following procedure {\em Steps 1--3}.

 {\em Step 1.} Choose a non-negative integer $p_0$ such that
\begin{equation} 
 \vPhi(x^{-1})\cap\wt{\vPhi}^\infty(\bs)\subset
 \vPhi(z_{\bs}\ft{(p_0)}).
\end{equation}

 {\em Step 2.} In the case where $xz_{\bs}\ft{(p_0)}=1$, put
 $\bs'\ft{(p)}:=\bs\ft{(p_0+p)}$ for each $p\in{\mathbb N}$.
 In the case where $xz_{\bs}\ft{(p_0)}\neq1$,
 choose a reduced expression
 $xz_{\bs}\ft{(p_0)}=\bs'\ft{(1)}\cdots\bs'\ft{(l_0)}$ with
 $l_0\in{\mathbb N}$, and put
 $\bs'\ft{(l_0+p)}:=\bs\ft{(p_0+p)}$ for each $p\in{\mathbb N}$.

 {\em Step 3.} Set $\bs':=(\bs'\ft{(p)})_{p\in{\mathbb N}}$.

 {\em(3)} If $x\in W$ and $\bs\in{\cl W}^\infty$, then
\[ \wt{\vPhi}^\infty(x,\bs)=\{\vPhi(x)\setminus(-\Omega)\}
 \amalg\{x\wt{\vPhi}^\infty(\bs)\setminus\Omega\}, \]
 where $\Omega:=x\wt{\vPhi}^\infty(\bs)\cap\nre$.
 In particular, if $x\wt{\vPhi}^\infty(\bs)\subset\pre$ then
\[ \wt{\vPhi}^\infty(x,\bs)=
 \vPhi(x)\amalg{x\wt{\vPhi}^\infty(\bs)}. \]

 {\em(4)} Suppose that $\bs,\,\bs'\in{\cl W}^\infty$.
\begin{itemize}
\item[(i)] If $\wt{\vPhi}^\infty(\bs)=\wt{\vPhi}^\infty(\bs')$,
 then $\wt{\vPhi}^\infty(x,\bs)=\wt{\vPhi}^\infty(x,\bs')$
 for each $x\in W$,
\item[(ii)] If $\wt{\vPhi}^\infty(y,\bs)=\wt{\vPhi}^\infty(\bs')$,
 then $\wt{\vPhi}^\infty(xy,\bs)=\wt{\vPhi}^\infty(x,\bs')$
  for each $(x,y)\in W^2$.
\end{itemize}
\end{lem}

\begin{proof}
 (1) Suppose that $\be\in\wt{\vPhi}^\infty(1,\bs)$.
 Since $z_{\bs}\ft{(p_0)}^{-1}(\be)<0$
 for some $p_0\in{\mathbb N}$, we have
 $\be\in\wt{\vPhi}^\infty(z_{\bs}\ft{(p_0)})\subset
 \wt{\vPhi}^\infty(\bs)$. Thus we get
 $\wt{\vPhi}^\infty(1,\bs)\subset\wt{\vPhi}^\infty(\bs)$.
 On the other hand, for each $p<q$, we have
\[ z_{\bs}\ft{(q)}^{-1}(\phi_{\bs}\ft{(p)})
 =-s_{\bs(q)}\cdots s_{\bs(p+1)}(\al_{\bs(p)})<0, \]
 and hence $\phi_{\bs}\ft{(p)}\in\wt{\vPhi}^\infty(1,\bs)$
 for each $p\in{\mathbb N}$. Thus we get
 $\wt{\vPhi}^\infty(\bs)\subset\wt{\vPhi}^\infty(1,\bs)$
 by Proposition 2.10(2).

 (2) Let $\bs'$ be an element of $S^{\mathbb N}$
 constructed as in (Step1)--(Step3).
 By the construction, we have
 $xz_{\bs}\ft{(p_0)}=z_{\bs'}\ft{(l_0)}$
 for some unique $l_0\in{\mathbb Z}_{\geq0}$.
 Since $\bs\ft{(p_0+p)}=\bs'\ft{(l_0+p)}$
 for each $p\in{\mathbb N}$, we have
\begin{align}
 xz_{\bs}\ft{(p_0+p)}&=z_{\bs'}\ft{(l_0+p)}, \\ 
 x\phi_{\bs}\ft{(p_0+p)}&=\phi_{\bs'}\ft{(l_0+p)}. 
\end{align}
 By the condition (2.1) and the equality (2.3),
 we have $\phi_{\bs'}\ft{(l_0+p)}>0$ for each $p\in{\mathbb N}$
 since $\phi_{\bs}\ft{(p_0+p)}\notin\vPhi(z_{\bs}\ft{(p_0)})$.
 In addition, by Theorem 2.2 we have $\phi_{\bs'}\ft{(p)}>0$
 for each $1\leq p\leq l_0$.
 Thus we get $\bs'\in{\cl W}^\infty$ by Proposition 2.10(1).
 Moreover, by (1) and the equality (2.2), we get
 $\wt{\vPhi}^\infty(x,\bs)=\wt{\vPhi}^\infty(\bs')$.

 (3) Since
 $-\Ome=\vPhi(x)\cap(-x\wt{\vPhi}^\infty(\bs))$, we have
\[ \vPhi(x)\setminus(-\Ome)=\{\, \be\in\pre \;|\; \be\in\vPhi(x),\;
 -x^{-1}(\be)\in\pre\setminus\wt{\vPhi}^\infty(\bs) \,\}. \]
 On the other hand, since
 $x\wt{\vPhi}^\infty(\bs)\setminus\Ome=
 x\wt{\vPhi}^\infty(\bs)\cap\pre$, we have
\[ x\wt{\vPhi}^\infty(\bs)\setminus\Ome=
 \{\,\be\in\pre\;|\;\be\notin\vPhi(x),\;
 x^{-1}(\be)\in\wt{\vPhi}^\infty(\bs)\,\}. \]
 Therefore, by (1) we get
 $\wt{\vPhi}^\infty(x,\bs)=\{\vPhi(x)\setminus(-\Ome)\}\amalg
 \{x\wt{\vPhi}^\infty(\bs)\setminus\Ome\}$.

 (4)(i) This is straightforward from (3).

 (4)(ii) By the argument in the proof of (2),
 there exist an element $\til{\bs}\in{\cl W}^\infty$ and
 $(p_0,l_0)\in({\mathbb Z}_{\geq0})^2$ satisfying
 $yz_{\bs}\ft{(p_0+p)}=z_{\til{\bs}}\ft{(l_0+p)}$
 for all $p\in{\mathbb N}$.
 Then we have
 $\wt{\vPhi}^\infty(\til{\bs})=
 \wt{\vPhi}^\infty(y,\bs)=\wt{\vPhi}^\infty(\bs')$.
 Hence, by (4)(i) we have
 $\wt{\vPhi}^\infty(x,\til{\bs})=\wt{\vPhi}^\infty(x,\bs')$.
 Moreover, since
 $xyz_{\bs}\ft{(p_0+p)}=xz_{\til{\bs}}\ft{(l_0+p)}$
 for all $p\in{\mathbb N}$, we have
 $\wt{\vPhi}^\infty(xy,\bs)=\wt{\vPhi}^\infty(x,\til{\bs})$.
 Thus we get
 $\wt{\vPhi}^\infty(xy,\bs)=
 \wt{\vPhi}^\infty(x,\til{\bs})=\wt{\vPhi}^\infty(x,\bs')$.
\end{proof}

\begin{defn} 
 Thanks to Proposition 2.10(4) and Lemma 2.13(1)(2)(4),
 we have a left action of $W$ on $W^\infty$ such that
 $x.[\bs]=[\bs']$ if $x\in W$ and $\bs,\,\bs'\in{\cl W}^\infty$
 satisfy $\wt{\vPhi}^\infty(x,\bs)=\wt{\vPhi}^\infty(\bs')$.
\end{defn}

\begin{prop} 
 If $x\in W$ and $\bs\in{\cl W}^\infty$, then
\[ \vPhi^\infty(x.[\bs])=\{\vPhi(x)\setminus(-\Omega)\}\amalg
 \{x\vPhi^\infty([\bs])\setminus\Omega\}, \]
 where $\Omega:=x\vPhi^\infty([\bs])\cap\nre$.
 In particular, if $x\vPhi^\infty([\bs])\subset\pre$ then
\[ \vPhi^\infty(x.[\bs])=\vPhi(x)\amalg{x\vPhi^\infty([\bs])}. \]
\end{prop}

\begin{proof}
 This follows from Lemma 2.13(3).
\end{proof}


%% file: sec3.tex
\section{Notation for the untwisted affine cases}
 In this section, we introduce notation for the untwisted affine cases.
 Let $l$ be a positive integer, and put
 ${\bd I}=\{0,1,\cdots,l\}$ and ${\SI}=\{1,\cdots,l\}$.
 Let ${\mathrm A}=(a_{ij})_{i,j\in{\bd I}}$ be a generalized
 Cartan matrix of affine type $X_l^{(1)}$ such that
 $(a_{ij})_{i,j\in{\SI}}$ is the Cartan matrix of type $X_l$,
 where $X=A,B,\cdots,G$.
 Let $({\fr h},\Pi,\Pi^\vee)$
 be a minimal realization of ${\mathrm A}$ over ${\mathbb R}$,
 that is, a triplet consisting of
 a $(l+2)$-dimensional real vector space ${\fr h}$
 and linearly independent subsets
 $\Pi=\{\al_i\,|\,i\in{\bd I}\}\subset{\fr h}^*$,
 $\Pi^\vee=\{\al_i^\vee\,|\,i\in{\bd I}\}\subset{\fr h}$
 satisfying $\la{\,\al_i^\vee\,,\,\al_j\,}\ra=a_{ij}$
 for each $i,\,j\in{\bd I}$.
 Let ${\fr g}$ be the affine Kac-Moody Lie algebra associated with
 $({\fr h},\Pi,\Pi^\vee)$, $\De\subset{\fr h}^*\setminus\{0\}$
 the root system of ${\fr g}$, $\re$ (resp. $\im$)
 the set of all real (resp. imaginary) roots, and
 $W=\la{\,s_i\,|\,i\in{\bd I}\,}\ra\subset\mathrm{GL}({\fr h}^*)$
 the Weyl group of ${\fr g}$, where $s_i$ is the reflection
 associated with $\al_i$. Furthermore, let $\De_+$ (resp. $\De_-$) be
 the set of all positive (resp. negative) roots relative to $\Pi$, and
 $\mathrm{ht}\colon\De_+\to{\mathbb N}$ the height function on $\De_+$.

 Set ${\sPi}:=\{\al_i\,|\,i\in{\SI}\}$,
 ${\st{\circ}{\fr h}}{}^*:=\oplus_{i\in{\SI}}{\mathbb R}\al_i$,
 ${\sW}:=\la{\,s_i\,|\,i\in{\SI}\,}\ra$,
 ${\sDe}:={\sW}(\st{\circ}{\Pi})$,
 and ${\sDe}_\pm:={\sDe}\cap\De_\pm$.
 Note that $\sDe$ is a root system of type $X_l$ with $\sPi$
 a root basis and $\sW$ is the Weyl group of $\sDe$.
 By the assumption on ${\mathrm A}$, we have
\[ \re=\{\,m\de+\vep\;|\;m\in{\mathbb Z},\,\vep\in{\sDe}\,\},
 \quad\im=\{\,m\de\;|\;m\in{\mathbb Z}\setminus\{0\}\,\}, \]
 where $\de=\al_0+\theta$ with $\theta$ the highest root of $\sDe$.

 Let $(\cdot\,|\,\cdot)$ be the standard symmetric bilinear
 form on ${\fr h}^*$, scaled so that $(\al\,|\,\al)=2$
 for all long roots $\al$ of ${\sDe}$.
 Note that $(\de\,|\,\al_i)=0$ for all $i\in{\bd I}$ and
 that the restriction of the form $(\cdot\,|\,\cdot)$
 to ${\st{\circ}{\fr h}}{}^*$ is positive-definite.
 For each $\lam\in{\fr h}^*$, we denote by $\ol{\lam}$
 the image of $\lam$ by the orthogonal projection onto
 ${\st{\circ}{\fr h}}{}^*$.
 Note that each $\be\in\De$ can be uniquely written as
 $m\de+\ol{\be}$ with $m\in{\mathbb Z}$ and
 $\ol{\be}\in{\sDe}\amalg\{0\}$.

 For each $\al\in\re$, we denote by $s_\al$
 the reflection with respect to $\al$.
 For each $\lam\in{\fr h}^*$, we define an element
 $t_{\lam}\in\mathrm{GL}({\fr h}^*)$ by setting
\[ t_{\lam}(\mu)=\mu+(\mu\,|\,\de)\lam
 -\{(\mu\,|\,\lam)+\textstyle{\frac{1}{2}}
 (\lam\,|\,\lam)(\mu\,|\,\de)\}\de \]
 for each $\mu\in{\fr h}^*$.
 We have $t_{\lam}(\mu)=\mu-(\mu\,|\,\lam)\de$
 for each $\mu\in{\fr h}_0^*$, where
 ${\fr h}_0^*:=\oplus_{i\in{\bd I}}{\mathbb R}\al_i$.

\begin{lem}[\cite{vK}] 
 Set $\check{\al}_i=\frac{2\al_i}{(\al_i\,|\,\al_i)}$
 for each $i\in{\SI}$, and set
 ${\sQ}{}^\vee:=\oplus_{i\in{\SI}}{\mathbb Z}\check{\al}_i$
 and $T=\{\,t_{\lam}\,|\,\lam\in{\sQ}{}^\vee\,\}$.
 Then $T$ is a normal subgroup of $W$ such that $W={\sW}\ltimes T$.
\end{lem}

 For each $x\in W$, we denote by $\ol{x}$ the unique element of
 $\sW$ such that $x\in\ol{x}T$.
 The mapping $\ol{\,\cdot\,}\colon W\to{\sW},\,x\mapsto\ol{x}$,
 is a group homomorphism, which satisfies that
 $\ol{x(\lam)}=\ol{x}(\ol{\lam})$ and $\ol{s_{\al}}=s_{\ol{\al}}$
 for each $x\in W$, $\lam\in{\fr h}_0^*$, and $\al\in\re$.

%% file: sec4.tex
\section{Preliminary results for classical root systems}
 In this section, we give preliminary results for
 classical root systems.
 We use the notation introduced in Section 3.
 For each subset ${\bd J}\subset{\SI}$, we set
\begin{align*}
 {\sPi}_{\tiJ}&:=\{\,\al_j\;|\;j\in{\bd J}\,\}\subset{\sPi},\qquad
 {\sW}_{\tiJ}:=\la{\,s_j\;|\;j\in{\bd J}\,}\ra\subset{\sW}, \\
 {\sDe}_{\tiJ}&:={\sW}_{\tiJ}({\sPi}_{\tiJ})\subset{\sDe},
 \quad\mbox{ and }\quad{\sDe}_{\tiJ\pm}={\sDe}_{\tiJ}\cap{\sDe}_\pm.
\end{align*}
 Note that ${\sDe}_{\tiJ}$ is a root system with ${\sPi}_{\tiJ}$
 a root basis and ${\sW}_{\tiJ}$ the Weyl group if
 ${\bd J}\neq\emptyset$.
 For each ${\bd K}\subset{\bd J}$, we denote by
 ${\sW}{}^{\tiK}_{\tiJ}$ the minimal coset representatives of
 the set ${\sW}_{\tiJ}/{\sW}_{\tiK}$ of all right cosets.
 If ${\bd J}={\SI}$ we denote it simply by ${\sW}{}^{\tiK}$.
 Note that each element $w\in{\sW}_{\tiJ}$ can be uniquely written
 as $w^{\tiK}w_{\tiK}$ with $w^{\tiK}\in{\sW}{}^{\tiK}_{\tiJ}$ and
 $w_{\tiK}\in{\sW}_{\tiK}$, where $w^{\tiK}$ is a unique element
 of the smallest length in the right coset $w{\sW}_{\tiK}$.
 Moreover, we have
\[ {\sW}{}^{\tiK}_{\tiJ}=\{\,w\in{\sW}_{\tiJ}\;|\;w(\al_j)>0
 \;\text{for all}\; j\in{\bd K}\,\}, \]
 and
 ${\sW}{}^{\tiK}_{\tiJ}{\sW}{}^{\tiL}_{\tiK}={\sW}{}^{\tiL}_{\tiJ}$
 if ${\bd L}\subset{\bd K}\subset{\bd J}$. In addition, we set
\[ {\sDe}{}^{\tiK}_{\tiJ}:={\sDe}_{\tiJ}\setminus{\sDe}_{\tiK},\qquad
 {\sDe}{}^{\tiK}_{\tiJ\pm}:={\sDe}{}^{\tiK}_{\tiJ}\cap{\sDe}_\pm. \]

\begin{lem} 
 {\em(1)} The following equality holds\/{\em:}
\[ {\sDe}{}^{\tiK}_{\tiJ+}=
 \{\,\vep=\sum_{j\in{\bd J}}m_j\al_j\in{\sDe}_{\tiJ+}\,
 (m_j\in{\mathbb Z}_{\geq0})\;|\;m_j>0\;\;\text{for some}\;\;
 j\in{\bd J}\setminus{\bd K} \,\}. \]

 {\em(2)} We have
 ${\sDe}{}^{\tiK}_{\tiJ\pm}\dotplus{\sDe}{}^{\tiK}_{\tiJ\pm}
 \subset{\sDe}{}^{\tiK}_{\tiJ\pm}$ and
 ${\sDe}{}^{\tiK}_{\tiJ\pm}\dotplus{\sDe}_{\tiK}
 \subset{\sDe}{}^{\tiK}_{\tiJ\pm}$.

 {\em(3)} For each $v\in{\sW}_{\tiK}$, we have
 $v{\sDe}{}^{\tiK}_{\tiJ\pm}={\sDe}{}^{\tiK}_{\tiJ\pm}$.

 {\em(4)} Let ${\bd K}_1$ and ${\bd K}_2$ be subsets of ${\bd J}$,
 and let $w_1$ and $w_2$ be elements of ${\sW}_{\tiJ}$.
 Then the following two conditions are equivalent{\em:}
\[ \text{\em(i)}\;\; w_1{\sDe}{}^{\tiK_1}_{\tiJ\pm}\subset
 w_2{\sDe}{}^{\tiK_2}_{\tiJ\pm};\qquad
 \text{\em(ii)}\;\;
 {\bd K}_1\supset{\bd K}_2,\; w_1\in w_2{\sW}_{\tiK_1}. \]
\end{lem}

\begin{proof}
 (1) This is straightforward from the definition.

 (2) This follows immediately from (1).

 (3) Let $\vep$ be an element of ${\sDe}{}^{\tiK}_{\tiJ+}$, and write
 $\vep=\sum_{j\in{\bd J}}m_j\al_j$ with $m_j\in{\mathbb Z}_{\geq0}$
 for all $j\in{\bd J}$ and $m_{j_\ast}>0$ for some
 $j_\ast\in{\bd J}\setminus{\bd K}$.
 Since $v(\al_j)\in\al_j+\sum_{k\in{\bd K}}{\mathbb Z}\al_k$
 for each $j\in{\bd J}\setminus{\bd K}$, we have
 $v(\vep)=\sum_{j\in{\bd J}\setminus{\bd K}}m_j\al_j+
 \sum_{k\in{\bd K}}m'_k\al_k\in{\sDe}$
 with $m'_k\in{\mathbb Z}$, which implies that
 $v(\vep)\in{\sDe}{}^{\tiK}_{\tiJ+}$ since $m_{j_\ast}>0$.
 Thus $v{\sDe}{}^{\tiK}_{\tiJ+}\subset{\sDe}{}^{\tiK}_{\tiJ+}$
 for each $v\in{\sW}_{\tiK}$, and hence
 $v{\sDe}{}^{\tiK}_{\tiJ+}={\sDe}{}^{\tiK}_{\tiJ+}$.

 (4) Suppose that ${\bd K}_1\supset{\bd K}_2$ and
 $w_1=w_2v$ with $v\in{\sW}_{\tiK_1}$.
 Then, by (3) we have
 $w_1{\sDe}{}^{\tiK_1}_{\tiJ+}=w_2{\sDe}{}^{\tiK_1}_{\tiJ+}
 \subset w_2{\sDe}{}^{\tiK_2}_{\tiJ+}$.
 Conversely, suppose that
 $w{\sDe}{}^{\tiK_1}_{\tiJ+}\subset{\sDe}{}^{\tiK_2}_{\tiJ+}$
 with $w=w_2^{-1}w_1$. Then we have
 $w^{\tiK_1}{\sDe}{}^{\tiK_1}_{\tiJ+}\subset
 {\sDe}{}^{\tiK_2}_{\tiJ+}$ by (3), and hence
 $w^{\tiK_1}(\al_j)>0$ for all $j\in{\bd J}\setminus{\bd K}_1$
 since
 ${\sPi}_{\ti{{\bd J}\setminus{\bd K}_1}}\subset
 {\sDe}{}^{\tiK_1}_{\tiJ+}$.
 Moreover, $w^{\tiK_1}(\al_k)>0$ for all $k\in{\bd K}_1$
 since $w^{\tiK_1}\in{\sW}{}^{\tiK_1}_{\tiJ}$.
 Thus $w^{\tiK_1}(\al_j)>0$ for all $j\in{\bd J}$, and hence
 $w^{\tiK_1}=1$ and $w=w_{\tiK_1}\in{\sW}_{\tiK_1}$.
 Therefore
 ${\sDe}{}^{\tiK_1}_{\tiJ+}=w{\sDe}{}^{\tiK_1}_{\tiJ+}
 \subset{\sDe}{}^{\tiK_2}_{\tiJ+}$,
 which implies that ${\bd K}_1\supset{\bd K}_2$.
\end{proof}

\begin{defn} 
 Let ${\bd J}$ be a non-empty subset of ${\SI}$.
\begin{itemize}
\item[(1)] A subset $P\subset{\sDe}$ is called
 a {\it closed\/} set if it satisfies the condition that
 if $\vep,\,\eta\in P,\;\vep+\eta\in\sDe$ then
 $\vep+\eta\in P$ (cf. \cite[\S1.7]{nB}).
 We call a subset $P\subset{\sDe}_{\tiJ}$
 a {\it coclosed set in\/} ${\sDe}_{\tiJ}$ if
 ${\sDe}_{\tiJ}\setminus P$ is a closed set, and call
 a subset $P\subset{\sDe}_{\tiJ}$ a {\it biclosed set in\/}
 ${\sDe}_{\tiJ}$ if both $P$ and ${\sDe}_{\tiJ}\setminus P$ are closed sets.
\item[(2)] We call a subset $P\subset{\sDe}_{\tiJ}$
 a {\it parabolic set in\/} ${\sDe}_{\tiJ}$ if
 $P$ is a closed set such that $P\cup(-P)={\sDe}_{\tiJ}$
 (cf. \cite{nB}).
\item[(3)] A subset $P\subset{\sDe}$ is called
 a {\it symmetric\/} set if $P=-P$ (cf. \cite{nB}).
\item[(4)] We call a subset $P\subset{\sDe}$
 a {\it pointed\/} set if $P\cap(-P)=\emptyset$.
\end{itemize}
\end{defn}

\begin{prop}[\cite{nB}] 
 The following three conditions are equivalent\/{\em:}
\begin{itemize}
\item[(i)] $P$ is a parabolic set in ${\sDe}_{\tiJ}${\em;}
\item[(ii)] $P$ is a closed subset of ${\sDe}_{\tiJ}$ such that
 $P\supset w{\sDe}_{\tiJ+}$ for some $w\in{\sW}_{\tiJ}${\em;}
\item[(iii)] $P=w({\sDe}_{\tiJ+}\amalg{\sDe}_{\tiK-})$
 for some ${\bd K}\subset{\bd J}$ and $w\in{\sW}_{\tiJ}$.
\end{itemize}
\end{prop}

\begin{prop}[\cite{nB}] 
 If $P$ is a pointed closed subset of ${\sDe}_{\tiJ}$, then
 there exists an element $w\in{\sW}_{\tiJ}$ such that
 $wP\subset{\sDe}_{\tiJ-}$.
\end{prop}

\begin{prop} 
 Let $P$ be a subset of ${\sDe}$.
 Then there exist a unique symmetric subset $P_s\subset P$ and
 a unique pointed subset $P_p\subset P$ such that
 $P=P_p\amalg P_s$. Moreover, if $P$ is closed then
 both $P_s$ and $P_p$ are closed sets satisfying
\begin{equation} 
 P_p\dotplus P_s\subset P_p.
\end{equation}
\end{prop}

\begin{proof}
 Suppose that there exist a symmetric subset $P_s\subset P$
 and a pointed subset $P_p\subset P$ such that $P=P_p\amalg P_s$.
 Then we have
\begin{align} 
 P_s&=\{\,\vep\in P\;|\;-\vep\in P\,\}, \\
 P_p&=\{\,\vep\in P\;|\;-\vep\in{\sDe}\setminus P\,\}.
\end{align}
 This proves the uniqueness of the decomposition.
 On the other hand, it is easy to see that
 the above subsets give the desired decomposition of $P$.

 In addition, we suppose that $P$ is closed.
 Let $\vep$ and $\eta$ be elements of $P_s$ such that
 $\vep+\eta\in{\sDe}$. Then we have
 $\vep+\eta\in P$ and $-\vep,\,-\eta\in P$.
 Thus we get $-(\vep+\eta)\in P$, and hence $\vep+\eta\in P_s$.
 Therefore $P_s$ is closed.

 We next prove (4.1).
 Suppose that $\vep+\eta\in P_s$ for some
 $\vep\in P_p$ and $\eta\in P_s$.
 Then $\vep=(\vep+\eta)+(-\eta)\in P_s$,
 since $P_s$ is closed and $-\eta\in P_s$.
 This is a contradiction. Hence, (4.1) is valid.

 Suppose that $\vep+\eta\in P_s$ for some $\vep,\,\eta\in P_p$.
 Then, since $-\vep-\eta\in P_s$ we have
 $-\vep=\eta+(-\vep-\eta)\in P_p$ by (4.1).
 This contradicts to $ P_p\cap(-P_p)=\emptyset$.
 Thus we get $\vep+\eta\in P_p$ for each $\vep,\,\eta\in P_p$
 satisfying $\vep+\eta\in{\sDe}$.
 Therefore $P_p$ is closed.
\end{proof}

\begin{prop} 
 The following four conditions are equivalent\/{\em:}
\begin{itemize}
\item[(i)] $P$ is a pointed biclosed set in ${\sDe}_{\tiJ}${\em;}
\item[(ii)] $P$ is a pointed coclosed set in ${\sDe}_{\tiJ}${\em;}
\item[(iii)] $P$ is a subset of ${\sDe}_{\tiJ}$ such that
 ${\sDe}_{\tiJ}\setminus P$ is a parabolic set in
 ${\sDe}_{\tiJ}${\em;}
\item[(iv)] $P=u{\sDe}{}^{\tiK}_{\tiJ-}$ for some unique
 ${\bd K}\subset{\bd J}$ and unique $u\in{\sW}{}^{\tiK}_{\tiJ}$.
\end{itemize}
\end{prop}

\begin{proof}
 (i)$\Rightarrow$(ii) It is Clear.

 (ii)$\Rightarrow$(iii)
 It is clear that $P\subset{\sDe}_{\tiJ}$ and
 ${\sDe}_{\tiJ}\setminus P$ is a closed set.
 By Proposition 4.5, we have
\begin{equation} 
 {\sDe}_{\tiJ}=P_p\amalg P_s\amalg
 ({\sDe}_{\tiJ}\setminus P)_p\amalg({\sDe}_{\tiJ}\setminus P)_s,
\end{equation}
 where $P_s$ (resp. $({\sDe}_{\tiJ}\setminus P)_s$) is
 the symmetric part of $P$ (resp. ${\sDe}_{\tiJ}\setminus P$) and
 $P_p$ (resp. $({\sDe}_{\tiJ}\setminus P)_p$) is the pointed part
 of $P$ (resp. ${\sDe}_{\tiJ}\setminus P$).
 Then we have
\begin{equation} 
 -P_p=({\sDe}_{\tiJ}\setminus P)_p.
\end{equation}
 Indeed, if $\vep\in P_p$ then we have
 $-\vep\in{\sDe}_{\tiJ}\setminus P$ and $-(-\vep)\in P$ by (4.2),
 and hence $-\vep\in({\sDe}_{\tiJ}\setminus P)_p$ by (4.3).
 Thus $-P_p\subset({\sDe}_{\tiJ}\setminus P)_p$.
 Similarly we have $-({\sDe}_{\tiJ}\setminus P)_p\subset P_p$.

 By (4.5), we have
 $-({\sDe}_{\tiJ}\setminus P)=
 P_p\amalg({\sDe}_{\tiJ}\setminus P)_s$.
 Moreover, we have $P_s=\emptyset$ since $P$ is pointed.
 Thus we get
 ${\sDe}_{\tiJ}=
 -({\sDe}_{\tiJ}\setminus P)\cup({\sDe}_{\tiJ}\setminus P)$
 by (4.4), and hence ${\sDe}_{\tiJ}\setminus P$
 is a parabolic set in ${\sDe}_{\tiJ}$.

 (iii)$\Rightarrow$(iv)
 By Proposition 4.3, there exist a subset
 ${\bd K}\subset{\bd J}$ and an element $w\in{\sW}_{\tiJ}$
 such that
 ${\sDe}_{\tiJ}\setminus P=w({\sDe}_{\tiJ+}\amalg{\sDe}_{\tiK-})$.
 Then $P=w{\sDe}{}^{\tiK}_{\tiJ-}$ since $P\subset{\sDe}_{\tiJ}$,
 and hence $P=w^{\tiK}{\sDe}{}^{\tiK}_{\tiJ-}$ by Lemma 4.1(3).
 The uniqueness follows from Lemma 4.1(4).

 (iv)$\Rightarrow$(i)
 It is clear that $u{\sDe}{}^{\tiK}_{\tiJ-}$ is pointed.
 By Lemma 4.1(2), we have
 $u{\sDe}{}^{\tiK}_{\tiJ-}\dotplus u{\sDe}{}^{\tiK}_{\tiJ-}
 \subset u{\sDe}{}^{\tiK}_{\tiJ-}$, and hence
 $u{\sDe}{}^{\tiK}_{\tiJ-}$ is closed.
 Moreover, by Lemma 4.1(2) we have
 $u{\sDe}{}^{\tiK}_{\tiJ+}\dotplus u{\sDe}{}^{\tiK}_{\tiJ+}
 \subset u{\sDe}{}^{\tiK}_{\tiJ+}$ and
 $u{\sDe}{}^{\tiK}_{\tiJ+}\dotplus u{\sDe}_{\tiK}
 \subset u{\sDe}{}^{\tiK}_{\tiJ+}$.
 In addition, $u{\sDe}_{\tiK}$ is closed.
 Thus ${\sDe}_{\tiJ}\setminus{u{\sDe}{}^{\tiK}_{\tiJ-}}$
 is closed, since
 ${\sDe}_{\tiJ}\setminus{u{\sDe}{}^{\tiK}_{\tiJ-}}=
 u{\sDe}{}^{\tiK}_{\tiJ+}\amalg u{\sDe}_{\tiK}$.
\end{proof}


%% file: sec5.tex
\section{The construction of biconvex sets}
 In this section, we give several methods of
 constructing biconvex sets for the root system
 of an arbitrary untwisted affine Lie algebra.

\begin{defn} 
 For each $\vep\in{\sDe}$ and $P\subset{\sDe}$, we define
 subsets $\la{\vep}\ra,\,\la{P}\ra\subset\pre$ by setting
\[ \la{\vep}\ra:=
 \{\,m\de+\vep\;|\;m\in{\mathbb Z}_{\geq0}\,\}\cap\pre,\qquad
 \la{P}\ra:=\coprod_{\vep\in{P}}\la{\vep}\ra. \]
\end{defn}

\begin{lem} 
 {\em(1)} Let $P$ be a subset of ${\sDe}$, and $x$ an element of $W$.
 Then
\[ \text{\em(i)}\;\;
 \la{P}\ra=\{\,\be\in\De_+\;|\;\ol{\be}\in P\,\};\quad
 \text{\em(ii)}\;\; \ol{x\la{P}\ra}\subset\ol{x}P;\quad
 \text{\em(iii)}\;\; x\la{P}\ra\doteq\la{\ol{x}P}\ra. \]

 {\em(2)} For subsets $P,\,P'\subset{\sDe}$,
 the following three conditions are equivalent\/{\em:}
\[ \text{\em(i)}\;\; P\subset{P'};\quad
 \text{\em(ii)}\;\; \la{P}\ra\subset\la{P'}\ra;\quad
 \text{\em(iii)}\;\; \la{P}\ra\,\dot{\subset}\,\la{P'}\ra. \]
\end{lem}

\begin{proof}
 (1) The (i) is straightforward from the definition.
 To prove (ii), suppose that $\be\in\la{P}\ra$.
 Then $\ol{\be}\in P$ by (i), hence
 $\ol{x(\be)}=\ol{x}(\ol{\be})\in\ol{x}P$. Thus (ii) is valid.
 We prove (iii). Write $x=t_{\lambda}\ol{x}$ with
 $\lambda\in{\sQ}{}^\vee$.
 Then we have
 $x(m\de+\vep)=(m-(\ol{x}\vep\,|\,\lambda))\de+\ol{x}\vep$
 for each $m\in{\mathbb Z}_{\geq0}$ and $\vep\in P$.
 Thus we get
 $x(m\de+\vep)\in\la{\ol{x}\vep}\ra$
 for all $m>(\ol{x}\vep\,|\,\lambda)$, and hence
 $x\la{\vep}\ra\,\dot{=}\,\la{\ol{x}\vep}\ra$.
 Thus (iii) is valid.

 (2) It is obvious.
\end{proof}

\begin{defn} 
 For each subset ${\bd J}\subset{\SI}$, we set
\begin{align*}
 \re_{\tiJ}&:=\la{{\sDe}_{\tiJ}}\ra,\qquad
 \De_{\tiJ}:=\re_{\tiJ}\amalg\im, \\
 \re_{\tiJ\pm}&:=\re_{\tiJ}\cap\De_\pm,\qquad
 \De_{\tiJ\pm}:={\De}_{\tiJ}\cap\De_\pm.
\end{align*}
 We denote by ${\fr B}_{\tiJ}$ the set of all finite biconvex sets
 in $\De_{\tiJ+}$, and by ${\fr B}_{\tiJ}^\infty$ the set of
 all infinite real biconvex sets in $\De_{\tiJ+}$.
 Note that $\De_{\tiJ+}$ is a convex set and that
 a subset $B\subset\De_{\tiJ+}$ is a convex set in $\De_{\tiJ+}$
 if and only if $B$ is a convex set (see Lemma 2.5(3)).
 For each non-empty subset ${\bd J}\subset{\SI}$, let
\[ {\sDe}_{\tiJ}=
 \coprod_{c=1}^{\ti{{\mathrm C}({\bd J})}}{\sDe}_{{\tiJ}_c} \]
 be the irreducible decomposition of ${\sDe}_{\tiJ}$ with
 $\ft{{\mathrm C}({\bd J})}$ a unique positive integer.
 Note that $\re_{\tiJ}=
 \coprod_{c=1}^{\ti{{\mathrm C}({\bd J})}}\re_{\tiJ_c}$.
 For each $c=1,\dots,\ft{{\mathrm C}({\bd J})}$,
 we denote by $\theta_{{\tiJ}_c}$ the highest root of
 ${\sDe}_{{\tiJ}_c}$ relative to the root basis
 ${\sPi}_{{\tiJ}_c}$, and set
\begin{align*}
 \Pi_{\tiJ_c}&:={\sPi}_{\tiJ_c}\amalg\{\de-\theta_{\tiJ_c}\}
 \;\;\text{for each}\;\;c=1,\dots,\ft{{\mathrm C}({\bd J})}, \\
 \Pi_{\tiJ}&:=\coprod_{c=1}^{\ti{{\mathrm C}({\bd J})}}
 {\Pi}_{\tiJ_c},\qquad S_{\tiJ}:=\{\,s_\al\;|\;\al\in\Pi_{\tiJ}\,\}.
\end{align*}
 We denote by $W_{\tiJ}$ the subgroup of $W$ generated by
 $S_{\tiJ}$, and by $V_{\tiJ}$ the subspace of ${\fr h}^*$
 spanned by $\Pi_{\tiJ}$. Note that
 $W_{\tiJ}=\prod_{c=1}^{\ti{{\mathrm C}({\bd J})}}W_{\tiJ_c}$
 (direct product) and that $\ol{y}\in{\sW}_{\tiJ}$
 for all $y\in W_{\tiJ}$.
 We also put $W_\emptyset:=\{1\}\subset W$.
 For each $s\in S_{\tiJ}$, we denote by $\al_s$
 the unique element of $\Pi_{\tiJ}$ such that $s=s_{\al_s}$.
\end{defn}

\begin{prop} 
 For each non-empty subset ${\bd J}\subset{\SI}$,
 the pair $(W_{\tiJ},S_{\tiJ})$ is a Coxeter system and
 the triplet $(V_{\tiJ},\De_{\tiJ},\Pi_{\tiJ})$ is
 a root system of $(W_{\tiJ},S_{\tiJ})$ with
 the properties {\em R(v)} and {\em R(vi)}.
\end{prop}

\begin{proof}
 Thanks to Theorem 2.2, it suffices to show that
 the triplet $(V_{\tiJ},\De_{\tiJ},\Pi_{\tiJ})$
 satisfies the six conditions R(i)--R(vi).
 The conditions R(i), R(iv), and
 the first equality in R(iii) are obvious.
 By definition, we have
\begin{equation*}
 \De_{\tiJ}=\De_{\tiJ-}\amalg\De_{\tiJ+},\qquad %
 \De_{\tiJ-}=-\De_{\tiJ+},\qquad \De_{\tiJ+}= %
 \bigcup_{c=1}^{\ti{{\mathrm C}({\bd J})}}\De_{\tiJ_c+}.
\end{equation*}
 Hence, to check the condition R(ii), it suffices to
 show that each element of $\De_{\tiJ_c+}$ can be written as
 $\sum_{\al\in\Pi_{\tiJ_c}}x_\al\al$ with
 $x_\al\geq0$ for all $\al\in\Pi_{\tiJ_c}$.
 We have
\begin{equation} 
 m\de-\theta_{\tiJ_c}=(\de-\theta_{\tiJ_c})+(m-1)\de
\end{equation}
 for each $m\in{\mathbb Z}_{\geq2}$, and
\begin{equation} 
 m\de-\vep=\{m\de-(\vep+\al_j)\}+\al_j
\end{equation}
 for each $m\in{\mathbb Z}_{\geq1}$ and
 $\vep\in{\sDe}_{{\tiJ_c}+}\setminus\{\theta_{{\tiJ}_c}\}$,
 where $j\in{\bd J}_c$ such that $\vep+\al_j\in{\sDe}_{{\tiJ_c}+}$.
 If $\be\in{\De}_{{\tiJ_c}+}\setminus{\Pi_{\tiJ_c}}$ satisfies
 $\ol{\be}\in{\sDe}_{{\tiJ_c}+}$, then we have either
\begin{equation} 
 \be=m\de+\ol{\be}\;\;\mbox{with $m\geq1$}\quad\mbox{or}\quad %
 \be=\ol{\be}=\vep+\eta\;\; %
 \mbox{with $\vep,\,\eta\in{\sDe}_{{\tiJ_c}+}$}.
\end{equation}
 In addition, we have
\begin{equation} 
 \de=(\de-\theta_{\tiJ_c})+\theta_{\tiJ_c}\quad\mbox{and}\quad %
 m\de=(m-1)\de+\de
\end{equation}
 for each $m\in{\mathbb Z}_{\geq2}$. By (5.1)--(5.4), we see that
 each element of $\De_{\tiJ_c+}\setminus\Pi_{\tiJ_c}$ can be
 written as $\be+\ga$ with $\be,\,\ga\in\De_{\tiJ_c+}$.
 Hence, by induction on values of elements of $\De_{\tiJ_c+}$
 by the height function $\mathrm{ht}\colon\De_+\to{\mathbb N}$,
 we see that each element of $\De_{\tiJ_c+}$ can be written as
 a ${\mathbb Z}_{\geq0}$-linear combination of $\Pi_{\tiJ_c}$.
 Thus R(ii) and R(v) are satisfied, and R(vi) is clear since
 $\mathrm{ht}|_{\De_{\tiJ+}}$ satisfies the required property in R(vi).
 Finally, we check the second equality in R(iii).
 Suppose that $\al_{s_0}\in\Pi_{{\tiJ}_c}$ with $s_0\in S_{\tiJ_c}$.
 Since $\mbox{${\mathbb N}\al_{s_0}$}\cap\De_{\tiJ_c+}=\{\al_{s_0}\}$,
 each element of $\De_{\tiJ_c+}\setminus\{\al_{s_0}\}$ can be
 written as $\sum_{s\in S_{\tiJ_c}}x_s\al_s$ with $x_s\geq0$
 for all $s\in S_{\tiJ_c}$ and $x_{s_1}>0$ for some $s_1\neq s_0$.
 This fact implies that
 $s_0(\De_{\tiJ_c+}\setminus\{\al_{s_0}\})= %
 \De_{\tiJ_c+}\setminus\{\al_{s_0}\}$.
 Thus the second equality in R(iii) is valid, since
 $s_0$ fixes pointwise $\De_{\tiJ+}\setminus\De_{\tiJ_c+}$.
\end{proof}

\begin{cor} 
 Let ${\bd J}$ be an arbitrary non-empty subset of ${\SI}$.

 {\em(1)} The assignment
 $y\,\mapsto\,\vPhi_{\tiJ}(y):=\vPhi(y)\cap\De_{\tiJ+}$ defines
 a bijective mapping from $W_{\tiJ}$ to ${\fr B}_{\tiJ}$.

 {\em(2)} Suppose that $y=s_1s_2\cdots s_n$ with $n\in{\mathbb N}$
 and $s_1,s_2,\dots,s_n\in S_{\tiJ}$ is a reduced expression of
 an element $y\in W_{\tiJ}\setminus\{1\}$. Then
 the following equality holds\/{\em:}
\[ \vPhi_{\tiJ}(y)=\{\,\al_{s_1},\,s_1(\al_{s_2}),
 \,\dots,\,s_1\cdots s_{n-1}(\al_{s_n})\,\}, \]
 where the elements of $\vPhi_{\tiJ}(y)$ displayed above are
 distinct from each other. In particular,
 $\sharp\vPhi_{\tiJ}(y)=\ell_{\tiJ}(y)$, where
 $\ell_{\tiJ}\colon W_{\tiJ}\to{\mathbb Z}_{\geq0}$ is
 the length function of $(W_{\tiJ},S_{\tiJ})$.
\end{cor}

\begin{proof}
 Since
 $\vPhi_{\tiJ}(y)=\{\,\be\in\De_{\tiJ+}\;|\;y^{-1}(\be)<0\,\}$,
 the part (1) follows from Theorem 2.6 and Proposition 5.4.
 The part (2) follows from Theorem 2.2 and Proposition 5.4.
\end{proof}

\noindent
{\it Remarks.}
 (1) A assertion similar to the part (1) of Corollary 5.5
 was stated by P.~Cellini and P.~Papi in the proof of
 Theorem 3.12 in \cite{CP} with an outline of the proof.
 However, it seems that the detailed proof is not given
 in the paper.

 (2) The action of $W_{\tiJ}$ on $V_{\tiJ}$ is faithful.
 Indeed, if $y\ft{|_{V_{\bd J}}}=id_{V_{\bd J}}$ for $y\in W_{\tiJ}$,
 then $\vPhi_{\tiJ}(y)=\emptyset$, and hence $y=1$ by the part (2) of
 Corollary 5.5. Therefore we may identify $W_{\tiJ}$ with the subgroup
 of $\mathrm{GL}(V_{\tiJ})$ generated by
 $s\ft{|_{V_{\bd J}}}$ with $s\in S_{\tiJ}$.

 (3) Set
 ${\sQ}{}_{\tiJ}^\vee:=\oplus_{j\in{\bd J}}{\mathbb Z}\check{\al}_j$
 and $T_{\tiJ}:=
 \{\,t_\lam\ft{|_{V_{\bd J}}}\,|\,\lam\in{\sQ}{}_{\tiJ}^\vee\,\}$.
 Then $W_{\tiJ}={\sW}_{\tiJ}\ltimes T_{\tiJ}$ (see Lemma 3.1).

\begin{defn} 
 For each $w\in{\sW}_{\tiJ}$ and ${\bd K}\subset{\bd J}$, we set
\begin{equation*}
 \De^{\tiK}_{\tiJ}\ft{(w,\pm)}:=\la{w{\sDe}{}^{\tiK}_{\tiJ\pm}}\ra.
\end{equation*}
 We denote it simply by $\De_{\tiJ}\ft{(w,\pm)}$ if
 ${\bd K}=\emptyset$, by $\De^{\tiK}\ft{(w,\pm)}$ if
 ${\bd J}={\SI}$, and by $\De\ft{(w,\pm)}$ if
 ${\bd K}=\emptyset$ and ${\bd J}={\SI}$.
\end{defn}

\begin{lem} 
 {\em(1)} The set $\De^{\tiK}_{\tiJ}\ft{(w,\pm)}$ is an infinite set
 if and only if ${\bd K}\subsetneq{\bd J}$.

 {\em(2)} For each $u\in{\sW}{}^{\tiK}_{\tiJ}$ and
 $v\in{\sW}_{\tiK}$, we have
\begin{align} 
 \De^{\tiK}_{\tiJ}\ft{(uv,\pm)}&=
 \De^{\tiK}_{\tiJ}\ft{(u,\pm)}, \\
 \De^{\tiK}_{\tiJ}\ft{(u,-)}&=
 \vPhi(u)\amalg u\De^{\tiK}_{\tiJ}\ft{(1,-)}, \\
 \De_{\tiJ+}^{re}&=\De^{\tiK}_{\tiJ}\ft{(u,-)}\amalg u\re_{\tiK+}
 \amalg\De^{\tiK}_{\tiJ}\ft{(u,+)}.
\end{align}
\end{lem}

\begin{proof}
 (1) This follows from the fact that
 ${\sDe}{}^{\tiK}_{\tiJ\pm}$ is not empty
 if and only if ${\bd K}\subsetneq{\bd J}$.

 (2) By Lemma 4.1(3), we have
 $uv{\sDe}{}^{\tiK}_{\tiJ\pm}=u{\sDe}{}^{\tiK}_{\tiJ\pm}$,
 which implies (5.5). By definition, we have
\begin{align*}
 \De^{\tiK}_{\tiJ}\ft{(u,-)}&=
 ({\sDe}_{\tiJ+}\cap u{\sDe}{}^{\tiK}_{\tiJ-})\amalg
 \{\,m\de+\vep\;|\;m\in\mbox{${\mathbb N}$},\,
 \vep\in u{\sDe}{}^{\tiK}_{\tiJ-}\,\} \\
 &=({\sDe}_{\tiJ+}\cap u{\sDe}{}^{\tiK}_{\tiJ-})\amalg
 u\De^{\tiK}_{\tiJ}\ft{(1,-)}.
\end{align*}
 Moreover, since $u\in{\sW}{}^{\tiK}_{\tiJ}$ we have
 $u{\sDe}_{\tiK-}\subset{\sDe}_-$, and hence
 $\vPhi(u)={\sDe}_{\tiJ+}\cap u{\sDe}{}^{\tiK}_{\tiJ-}$.
 Thus (5.6) is valid. By definition, we have
\[ \re_{\tiJ+}=\la{{\sDe}_{\tiJ}}\ra,\quad
 \De^{\tiK}_{\tiJ}\ft{(u,\pm)}=\la{u{\sDe}{}^{\tiK}_{\tiJ\pm}}\ra,
 \quad u\re_{\tiK+}=\la{u{\sDe}_{\tiK}}\ra. \]
 Thus (5.7) is valid, since
 ${\sDe}_{\tiJ}=u{\sDe}{}^{\tiK}_{\tiJ-}\amalg
 u{\sDe}_{\tiK}\amalg u{\sDe}{}^{\tiK}_{\tiJ+}$.
\end{proof}

\begin{prop} 
 Let $P$ be a subset of ${\sDe}$, and ${\bd J}$
 a non-empty subset of ${\SI}$.
\begin{itemize}
\item[(1)] If $P$ is a closed set, then
 $\la{P}\ra\amalg\pim$ is a convex set.
\item[(2)] If $P$ is a pointed closed set, then
 $\la{P}\ra$ is a real convex set.
\item[(3)] If $P$ is a pointed biclosed set in ${\sDe}_{\tiJ}$,
 then $\la{P}\ra$ is a real biconvex set in $\De_{\tiJ+}$.
\end{itemize}
\end{prop}

\begin{proof}
 (1) Suppose that $\be+\ga\in\De_+$ with
 $\be,\,\ga\in\la{P}\ra\amalg\pim$. Then
 $\ol{\be},\,\ol{\ga}\in P\amalg\{0\}$ and
 $\ol{\be}+\ol{\ga}\in{\sDe}\amalg\{0\}$.
 Since $P$ is closed, we have $\ol{\be}+\ol{\ga}\in P\amalg\{0\}$,
 and hence $\be+\ga\in\la{P}\ra\amalg\pim$ by (i) of Lemma 5.2(1).
 Thus $\la{P}\ra\amalg\pim$ is a convex set.

 (2) It is clear that $\la{P}\ra\subset\pre$. Suppose that
 $\be+\ga\in\De_+$ with $\be,\,\ga\in\la{P}\ra$.
 Then $\ol{\be},\,\ol{\ga}\in P$ and
 $\ol{\be}+\ol{\ga}\in{\sDe}\amalg\{0\}$.
 If $\ol{\be}+\ol{\ga}=0$ then $\ol{\be}=-\ol{\ga}\in P\cap(-P)$.
 This contradicts to $P\cap(-P)=\emptyset$.
 Thus we get $\ol{\be}+\ol{\ga}\in{\sDe}$.
 Since $P$ is closed, we have $\ol{\be}+\ol{\ga}\in P$,
 and hence $\be+\ga\in\la{P}\ra$ by (i) of Lemma 5.2(1).
 Therefore $\la{P}\ra$ is a real convex set.

 (3) From (2), it follows that $\la{P}\ra$ is a real convex set.
 Since $P$ is a biclosed set in ${\sDe}_{\tiJ}$, the set
 ${\sDe}_{\tiJ}\setminus P$ is a closed set.
 Thus, by (1), we see that $\De_{\tiJ+}\setminus\la{P}\ra$
 is a convex set, since
 $\De_{\tiJ+}\setminus\la{P}\ra=
 \la{{\sDe}_{\tiJ}\setminus P}\ra\amalg\pim$.
\end{proof}

\begin{cor} 
 Let ${\bd K}$ be a subset of ${\bd J}$, and
 $u$ an element of ${\sW}{}^{\tiK}_{\tiJ}$. Then
\[ \text{{\em(i)}\;\;$u\De_{\tiK+}$ is a convex set\/};\qquad
 \text{{\em(ii)}\;\;$\De^{\tiK}_{\tiJ}\ft{(u,\pm)}$
 is a real biconvex set in $\De_{\tiJ+}$}. \]
\end{cor}

\begin{proof}
 We have $u\De_{\tiK+}=\la{u{\sDe}_{\tiK}}\ra\amalg\pim$.
 Since $u{\sDe}_{\tiK}$ is a closed set,
 (i) follows from Proposition 5.8(1).
 From Proposition 4.6, it follows that
 $u{\sDe}{}^{\tiK}_{\tiJ\pm}$ is a pointed biclosed set
 in ${\sDe}_{\tiJ}$, hence (ii) follows from Proposition 5.8(3).
\end{proof}

\begin{lem} 
 For ${\bd K}\subset{\bd J}$ and $u\in{\sW}{}^{\tiK}_{\tiJ}$,
 we have
 $\De^{\tiK}_{\tiJ}\ft{(u,\pm)}\dotplus u\De_{\tiK+}
 \subset\De^{\tiK}_{\tiJ}\ft{(u,\pm)}$.
\end{lem}

\begin{proof}
 Suppose that $\be+\ga\in\De_+$ with
 $\be\in\De^{\tiK}_{\tiJ}\ft{(u,\pm)}$ and $\ga\in u\De_{\tiK+}$.
 Then we have
 $\ol{\be}\in u{\sDe}{}^{\tiK}_{\tiJ\pm}$,
 $\ol{\ga}\in u{\sDe}_{\tiK}\amalg\{0\}$, and
 $\ol{\be}+\ol{\ga}\in{\sDe}\amalg\{0\}$.
 Thus we get
 $\ol{\be}+\ol{\ga}\in u{\sDe}{}^{\tiK}_{\tiJ\pm}$
 by Lemma 4.1(2), and hence
 $\be+\ga\in\De^{\tiK}_{\tiJ}\ft{(u,\pm)}$
 by (i) of Lemma 5.2(1).
\end{proof}

\begin{prop} 
 Let ${\bd K}$ be a subset of ${\bd J}$, and
 $u$ an element of ${\sW}{}^{\tiK}_{\tiJ}$.
\begin{itemize}
\item[(1)] If $C$ is a convex set in $u\De_{\tiK+}$, then
 $C\amalg\De^{\tiK}_{\tiJ}\ft{(u,\pm)}$ is a convex set
 in $\De_{\tiJ+}$.
\item[(2)] If $C$ is a biconvex set in $u\De_{\tiK+}$, then
 $C\amalg\De^{\tiK}_{\tiJ}\ft{(u,\pm)}$ is a biconvex set
 in $\De_{\tiJ+}$.
\end{itemize}
\end{prop}

\begin{proof}
 (1) From (ii) of Corollary 5.9, it follows that
 $\De^{\tiK}_{\tiJ}\ft{(u,\pm)}$ is a convex set in $\De_{\tiJ+}$.
 Thus the assertion follows from Lemma 2.5(4) and Lemma 5.10.

 (2) By the equality (5.7), we have
\[ \De_{\tiJ+}\setminus\{C\amalg\De^{\tiK}_{\tiJ}\ft{(u,-)}\}=
 (u\De_{\tiK+}\setminus C)\amalg\De^{\tiK}_{\tiJ}\ft{(u,+)}. \]
 Since both $C$ and $u\De_{\tiK+}\setminus C$ are convex sets in
 $u\De_{\tiK+}$, we see that both
 $C\amalg\De^{\tiK}_{\tiJ}\ft{(u,-)}$ and
 $(u\De_{\tiK+}\setminus C)\amalg\De^{\tiK}_{\tiJ}\ft{(u,+)}$
 are convex sets in $\De_{\tiJ}$ by (1), hence
 $C\amalg\De^{\tiK}_{\tiJ}\ft{(u,-)}$
 is a biconvex set in $\De_{\tiJ+}$.
 To prove of the assertion for $C\amalg\De^{\tiJ}\ft{(u,+)}$,
 it suffices to exchange the sign $\De^{\tiJ}\ft{(u,-)}$ for
 $\De^{\tiJ}\ft{(u,+)}$.
\end{proof}


%% file: sec6.tex
\section{A Parametrization of infinite real biconvex sets}
 In this section, we give a Parametrization of the set
 ${\fr B}_{\tiJ}^\infty$ of all infinite real biconvex sets in
 $\De_{\tiJ+}$ for each non-empty subset ${\bd J}\subset{\SI}$.

\begin{lem} 
 {\em(1)} If $B$ is a real coconvex set in $\De_{\tiJ+}$,
 then for each $\vep\in{\sDe}_{\tiJ}$ we have either
 $\la{\vep}\ra\subset B$ or
 $\la{\vep}\ra\,\dot{\subset}\,\re_{\tiJ+}\setminus B$.

 {\em(2)} If $B$ is a real biconvex set in $\De_{\tiJ+}$
 and a subset $P\subset{\sDe}_{\tiJ}$ satisfies
 $\la{P}\ra\,\dot{\subset}\,B$, then we have
 $\la{P}\ra\subset B$ and $\la{-P}\ra\cap B=\emptyset$.
\end{lem}

\begin{proof}
 (1) Suppose that there exists $m\in{\mathbb Z}_{\geq0}$
 such that $m\de+\vep\in\re_{\tiJ+}\setminus B$.
 Since $B\subset\pre$ we have
 $\im_+\subset\De_{\tiJ+}\setminus B$.
 Thus we get $(m+l)\de+\vep\in\re_{\tiJ+}\setminus B$
 for all $l\in{\mathbb N}$ by the convexity of
 $\De_{\tiJ+}\setminus B$, and hence
 $\la{\vep}\ra\,\dot{\subset}\,\re_{\tiJ+}\setminus B$.

 (2) By (1), we have $\la{P}\ra\subset B$.
 Suppose that $\la{-P}\ra\cap B\neq\emptyset$.
 Then there exists an element $\vep\in P$ such that
 $m\de-\vep\in B$ for some $m\in{\mathbb Z}_{\geq0}$.
 Moreover, we have $\de+\vep\in B$
 since $\la{\vep}\ra\subset B$.
 By the convexity of $B$, we have
 $(m+1)\de=(m\de-\vep)+(\de+\vep)\in B$.
 This contradicts to $B\subset\re_{\tiJ+}$.
 Hence we have $\la{-P}\ra\cap B=\emptyset$.
\end{proof}

\begin{prop} 
 Let $B$ be a real convex set in $\De_{\tiJ+}$, and set
\[ \ol{B}:=\{\,\ol{\be}\;|\;\be\in B\,\},\quad
 P_{\ti{B}}:=\{\,\vep\in{\sDe}_{\tiJ}\;|\;
 \la{\vep}\ra\,\dot{\subset}\,B\,\}. \]
 Then both $\ol{B}$ and $P_{\ti{B}}$ are pointed closed subsets
 of ${\sDe}_{\tiJ}$ such that $P_{\ti{B}}\subset\ol{B}$. Moreover,
 if $B$ is a real biconvex set in $\De_{\tiJ+}$ then
 $P_{\ti{B}}$ is a pointed biclosed set in ${\sDe}_{\tiJ}$.
\end{prop}

\begin{proof}
 It is clear that $\ol{B},\,P_{\ti{B}}\subset{\sDe}_{\tiJ}$.
 Suppose that $\vep+\eta\in{\sDe}$ with $\vep,\,\eta\in\ol{B}$.
 By definition, there exist $\be,\,\ga\in B$ such that
 $\ol{\be}=\vep,\,\ol{\ga}=\eta$.
 By the convexity of $B$, we have $\be+\ga\in B$,
 and hence $\vep+\eta=\ol{\be+\ga}\in\ol{B}$.
 Thus $\ol{B}$ is a closed set.
 Suppose that $\vep+\eta\in{\sDe}$ with $\vep,\,\eta\in P_{\ti{B}}$.
 By definition, we have
 $\la{\vep}\ra,\,\la{\eta}\ra\,\dot{\subset}\,B$, and hence
 there exist $m,\,n\in{\mathbb Z}_{\geq0}$ such that
 $(m+k)\de+\vep\in B$ and $(n+k)\de+\eta\in B$ for all
 $k\in{\mathbb Z}_{\geq0}$. By the convexity of $B$, we have
 $(m+n+k)\de+\vep+\eta\in B$ for all $k\in{\mathbb Z}_{\geq0}$.
 Thus we get $\la{\vep+\eta}\ra\,\dot{\subset}\,B$,
 and hence $\vep+\eta\in P_{\ti{B}}$.
 Therefore $P_{\ti{B}}$ is a closed set.
 Suppose that $\vep\in\ol{B}\cap(-\ol{B})$.
 Then we have $\vep,\,-\vep\in\ol{B}$.
 Hence we may assume that $\vep\in\ol{B}\cap{\sDe}_{\tiJ+}$.
 Then there exist $m\in{\mathbb Z}_{\geq0}$ and $n\in{\mathbb N}$
 such that $m\de+\vep,\;n\de-\vep\in B$.
 By the convexity of $B$, we have
 $(m+n)\de=(m\de+\vep)+(n\de-\vep)\in B$.
 This contradicts to $B\subset\re_{\tiJ+}$.
 Thus we get $\ol{B}\cap(-\ol{B})=\emptyset$.
 Moreover, by definition, we have $P_{\ti{B}}\subset\ol{B}$,
 and hence $P_{\ti{B}}\cap(-P_{\ti{B}})=\emptyset$.

 Next we prove the second assertion. It suffices to show that
 $P_{\ti{B}}$ is a coclosed set in $\sDe_{\tiJ}$.
 By the definition of $P_{\ti{B}}$ and Lemma 6.1(1), we see that
\begin{equation} 
 P_{\ti{B}}=\{\,\vep\in{\sDe}_{\tiJ}\;|\;\la{\vep}\ra\subset B\,\},
 \quad{\sDe}_{\tiJ}\setminus P_{\ti{B}}=\{\,\vep\in{\sDe}_{\tiJ}
 \;|\;\la{\vep}\ra\,\dot{\subset}\,\re_{\tiJ+}\setminus B\,\}.
\end{equation}
 Suppose that $\vep+\eta\in{\sDe}$ with
 $\vep,\,\eta\in{\sDe}_{\tiJ}\setminus P_{\ti{B}}$.
 Then $\la{\vep}\ra,\,\la{\eta}\ra\,\dot{\subset}\,
 \re_{\tiJ+}\setminus B$ by (6.1).
 By the convexity of $\De_{\tiJ+}\setminus B$, we have
 $\la{\vep+\eta}\ra\,\dot{\subset}\,\re_{\tiJ+}\setminus B$,
 and hence $\vep+\eta\in{\sDe}_{\tiJ}\setminus P_{\ti{B}}$.
 Thus $P_{\ti{B}}$ is a coclosed set in $\sDe_{\tiJ}$.
\end{proof}

\begin{prop} 
 Let ${\bd J}$ be an arbitrary non-empty subset of ${\SI}$.

 {\em(1)} If $B$ is a real convex set in $\De_{\tiJ+}$,
 then there exists an element $w\in{\sW}_{\tiJ}$ such that
 $B\subset\De_{\tiJ}\ft{(w,-)}$.

 {\em(2)} The assignment $w\,\mapsto\,\De_{\tiJ}\ft{(w,-)}$
 defines a bijective mapping from ${\sW}_{\tiJ}$ to the set
 ${\fr M}$ of all maximal real convex sets in $\De_{\tiJ+}$
 {\em(}relative to the inclusion relation\/{\em)}.
 Moreover, ${\fr M}$ coincides with the set of
 all maximal real biconvex sets in $\De_{\tiJ+}$.
\end{prop}

\begin{proof}
 (1) From Proposition 6.2, it follows that $\ol{B}$
 is a pointed closed subset of ${\sDe}_{\tiJ}$.
 Hence there exists an element $w\in{\sW}_{\tiJ}$ such that
 $\ol{B}\subset w{\sDe}_{\tiJ-}$ by Proposition 4.4.
 Then we have
\[ B\subset\la{\,\ol{B}\,}\ra\subset
 \la{\,w{\sDe}_{\tiJ-}\,}\ra=\De_{\tiJ}\ft{(w,-)}. \]

 (2) From Corollary 5.9, it follows for each $w\in W$ that
 $\De_{\tiJ}\ft{(w,-)}$ is a real biconvex set in $\De_{\tiJ+}$.
 In particular, $\De_{\tiJ}\ft{(w,-)}$ is a real convex set
 in $\De_{\tiJ+}$ for each $w\in{\sW}$.
 To prove the maximality of $\De_{\tiJ}\ft{(w,-)}$, suppose that
 $\De_{\tiJ}\ft{(w,-)}\subset B$ for some real convex set $B$
 in $\De_{\tiJ+}$.
 By (1), there exists an element $w'\in{\sW}_{\tiJ}$ such that
 $B\subset\De_{\tiJ}\ft{(w',-)}$.
 Since $\De_{\tiJ}\ft{(w,-)}\subset\De_{\tiJ}\ft{(w',-)}$,
 it follows that $w{\sDe}_{\tiJ-}\subset w'{\sDe}_{\tiJ-}$,
 which implies that $w=w'$, and hence $\De_{\tiJ}\ft{(w,-)}=B$.
 Therefore $\De_{\tiJ}\ft{(w,-)}$ is a maximal real convex set
 in $\De_{\tiJ+}$.
 Moreover, by the above argument,
 the injectivity of the mapping is obvious.
 Finally, we prove the surjectivity of the mapping.
 Let $B$ be a maximal real convex set in $\De_{\tiJ+}$.
 By (1), there exists an element $w\in{\sW}_{\tiJ}$ such that
 $B\subset\De_{\tiJ}\ft{(w,-)}$.
 The maximality of $B$ implies that $B=\De_{\tiJ}\ft{(w,-)}$.
\end{proof}

\begin{prop} 
 Let ${\bd J}$ be an arbitrary non-empty subset of ${\SI}$,
 and $B$ a real biconvex set in $\De_{\tiJ+}$.
 Then there exist a unique subset ${\bd K}\subset{\bd J}$ and
 a unique element $u\in{\sW}{}^{\tiK}_{\tiJ}$ such that
 $\De^{\tiK}_{\tiJ}\ft{(u,-)}\subset B$ and
 $B\,\dot{\subset}\,\De^{\tiK}_{\tiJ}\ft{(u,-)}$.
 Moreover, $B$ is an infinite set if and only if
 ${\bd K}\subsetneq{\bd J}$.
\end{prop}

\begin{proof}
 From Proposition 6.2, it follows that $P_{\ti{B}}$
 is a pointed biclosed subset of ${\sDe}_{\tiJ}$.
 Hence, by Proposition 4.6,
 there exist a unique subset ${\bd K}\subset{\bd J}$
 and a unique element $u\in{\sW}{}^{\tiK}_{\tiJ}$ such that
 $P_{\ti{B}}=u{\sDe}{}^{\tiK}_{\tiJ-}$.
 By (6.1), we see that $\la{\vep}\ra\subset B$
 for each $\vep\in u{\sDe}{}^{\tiK}_{\tiJ-}$ and
 that $\la{\vep}\ra\,\dot{\subset}\,\re_{\tiJ+}\setminus B$ for each
 $\vep\in{\sDe}_{\tiJ}\setminus u{\sDe}{}^{\tiK}_{\tiJ-}$.
 Thus we get
 $\De^{\tiK}_{\tiJ}\ft{(u,-)}\subset B$ and
 $B\,\dot{\subset}\,\De^{\tiK}_{\tiJ}\ft{(u,-)}$.
 The second assertion follows from Lemma 5.7(1).
\end{proof}

\begin{defn} 
 For each non-empty subset ${\bd J}\subset{\SI}$, we set
\begin{align*}
 \wt{\bol{\cl P}}_{\tiJ}&:=\{\,(\scK,u,y)\;|\;
 \mbox{${\bd K}\subset{\bd J}$},\;
 u\in{\sW}{}^{\tiK}_{\tiJ},\;y\in W_{\tiK}\,\}, \\
 \bol{\cl P}_{\tiJ}&:=\{\,(\scK,u,y)\in\wt{\bol{\cl P}}_{\tiJ}\;|\;
 \mbox{${\bd K}\subsetneq{\bd J}$}\,\},
\end{align*}
 where $W_{\tiK}$ is the subgroup of $W$ defined in Definition 5.3.
 For each $(\scK,u,y)\in\wt{\bol{\cl P}}_{\tiJ}$, we define a subset
 $\nab_{\tiJ}\ft{(\tiK,u,y)}\subset\re_{\tiJ+}$ by setting
\[ \nab_{\tiJ}\ft{(\tiK,u,y)}:=
 \De^{\tiK}_{\tiJ}\ft{(u,-)}\amalg u\vPhi_{\tiK}(y). \]
 Note that $\nab_{\tiJ}\ft{(\tiK,u,y)}=\vPhi_{\tiJ}(y)$
 if ${\bd K}={\bd J}$ and that
 $\nab_{\tiJ}\ft{(\tiK,u,y)}=\De_{\tiJ}\ft{(u,-)}$
 if ${\bd K}=\emptyset$.
\end{defn}

\begin{lem} 
 {\em(1)} For each $(\scK,u,y)\in\wt{\bol{\cl P}}_{\tiJ}$,
 the following equality holds\/{\em:}
\begin{equation}
 \nab_{\tiJ}\ft{(\tiK,u,y)}=
 \vPhi(u)\amalg u\nab_{\tiJ}\ft{(\tiK,1,y)}.
\end{equation}
 Moreover, $\nab_{\tiJ}\ft{(\tiK,u,y)}$
 is an infinite set if and only if
 $(\scK,u,y)\in\bol{\cl P}_{\tiJ}$.

 {\em(2)} Let $(\scK_1,u_1,y_1)$ and $(\scK_2,u_2,y_2)$ be
 elements of $\wt{\bol{\cl P}}_{\tiJ}$.
 Then the following two conditions are equivalent\/{\em:}
\[ \text{\em(i)}\;\;
 \nab_{\tiJ}\ft{(\tiK_1,u_1,y_1)}\,\;\dot{\subset}\;\,
 \nab_{\tiJ}\ft{(\tiK_2,u_2,y_2)};\qquad
 \text{\em(ii)}\;\;
 {\bd K}_1\supset{\bd K}_2,\;u_1\in u_2{\sW}_{\tiK_1}. \]
\end{lem}

\begin{proof}
 (1) By the equality (5.6), we have
\begin{align*}
 \vPhi(u)\amalg u\nab_{\tiJ}\ft{(\tiK,1,y)}&=
 \vPhi(u)\amalg u\De^{\tiK}_{\tiJ}\ft{(1,-)}\amalg
 u\vPhi_{\tiK}(y) \\
 &=\De^{\tiK}_{\tiJ}\ft{(u,-)}\amalg u\vPhi_{\tiK}(y)
 =\nab_{\tiJ}\ft{(\tiK,u,y)}.
\end{align*}
 The second assertion follows from Lemma 5.7(1).

 (2) The assertion follows from Lemma 4.1(4) and Lemma 5.2(2),
 since (i) is equivalent to the condition:
 $\De^{\tiK_1}_{\tiJ}\ft{(u_1,-)}\,\dot{\subset}\,
 \De^{\tiK_2}_{\tiJ}\ft{(u_2,-)}$.
\end{proof}

\begin{thm} 
 The assignment $(\scK,u,y)\mapsto\nab_{\tiJ}\ft{(\tiK,u,y)}$
 defines a bijective mapping from $\wt{\bol{\cl P}}_{\tiJ}$ to
 ${\fr B}_{\tiJ}\amalg{\fr B}_{\tiJ}^\infty$, which maps
 $\bol{\cl P}_{\tiJ}$ onto ${\fr B}_{\tiJ}^\infty$.
\end{thm}

\begin{proof}
 For each $(\scK,u,y)\in\wt{\bol{\cl P}}_{\tiJ}$, we see that
 $u\vPhi_{\tiK}(y)$ is a biconvex set in $u\De_{\tiK+}$, and hence
 $\nab_{\tiJ}\ft{(\tiK,u,y)}$ is a real biconvex set in $\De_{\tiJ+}$
 by Proposition 5.11(2).
 Thus the mapping $\nab_{\tiJ}$ is well-defined.
 Moreover, by the second assertion in Lemma 6.6(1), we see that
 $\nab_{\tiJ}(\bol{\cl P}_{\tiJ})\subset{\fr B}_{\tiJ}^\infty$ and
 $\nab_{\tiJ}(\wt{\bol{\cl P}}_{\tiJ}\setminus\bol{\cl P}_{\tiJ})
 \subset{\fr B}_{\tiJ}$.
 To prove the injectivity, suppose that
 $\nab_{\tiJ}\ft{(\tiK_1,u_1,y_1)}=\nab_{\tiJ}\ft{(\tiK_2,u_2,y_2)}$.
 By Lemma 6.6(2), we have
 ${\bd K}_1={\bd K}_2$ and $u_1\in u_2{\sW}_{\tiK_1}$, and hence
 $u_1=u_2$ since $u_1,\,u_2\in{\sW}{}^{\tiK_1}_{\tiJ}$.
 Thus we get
 $\De^{\tiK_1}_{\tiJ}\ft{(u_1,-)}=\De^{\tiK_2}_{\tiJ}\ft{(u_2,-)}$
 and $\vPhi_{\tiK_1}(y_1)=\vPhi_{\tiK_1}(y_2)$.
 By Corollary 5.5(1), we get $y_1=y_2$ and
 $(\scK_1,u_1,y_1)=(\scK_2,u_2,y_2)$.
 Finally, we prove the surjectivity.
 Suppose that $B\in{\fr B}_{\tiJ}\amalg{\fr B}_{\tiJ}^\infty$.
 Then $B\subset\re_{\tiJ+}$. By Proposition 6.4, there exist a subset
 ${\bd K}\subset{\bd J}$ and an element $u\in{\sW}{}^{\tiK}_{\tiJ}$
 such that $\De^{\tiK}_{\tiJ}\ft{(u,-)}\subset B$ and
 $B\,\dot{\subset}\,\De^{\tiK}_{\tiJ}\ft{(u,-)}$.
 Then $B\cap u\re_{\tiK+}$ is a finite biconvex set in $u\De_{\tiK+}$,
 since $B\cap u\re_{\tiK+}=B\cap u\De_{\tiK+}$.
 By Corollary 5.5(1), there exists an element $y\in W_{\tiK}$
 such that $B\cap u\re_{\tiK+}=u\vPhi_{\tiK}(y)$.
 Moreover, we have
 $B\cap\De^{\tiK}_{\tiJ}\ft{(u,+)}=\emptyset$
 by Lemma 6.1(2). Thus we get
 $(\scK,u,y)\in\wt{\bol{\cl P}}_{\tiJ}$ and
 $B=\De^{\tiK}_{\tiJ}\ft{(u,-)}\amalg u\vPhi_{\tiK}(y)=
 \nab_{\tiJ}\ft{(\tiK,u,y)}$ by (5.7).
\end{proof}


%% file: sec7.tex
\section{Main theorem}
 In this section, we describe in detail relationships between
 the set ${\cl W}_{\tiJ}^\infty$ of all infinite reduced words
 of the Coxeter system $(W_{\tiJ},S_{\tiJ})$ and the set
 ${\fr B}_{\tiJ}^\infty$ of all infinite real biconvex sets
 in $\De_{\tiJ+}$ for each non-empty subset ${\bd J}\subset{\SI}$.
 Let $W_{\tiJ}^\infty$ be the quotient set of
 ${\cl W}_{\tiJ}^\infty$ obtained by applying Definition 2.9 to
 the Coxeter system $(W_{\tiJ},S_{\tiJ})$, and
 $\vPhi_{\tiJ}^\infty\colon W_{\tiJ}^\infty\to{\fr B}_{\tiJ}^\infty$
 the injective mapping obtained by applying Definition 2.11 to
 the root system $(V_{\tiJ},\De_{\tiJ},\Pi_{\tiJ})$ of
 the Coxeter system $(W_{\tiJ},S_{\tiJ})$.

\begin{prop}[\cite{jB}] 
 Let ${\bd K}$ be a proper subset of ${\bd J}$, and
 $\lam$ an element of the lattice ${\sQ}{}_{\tiJ}^\vee$
 {\em(}see the remarks below {\em Corollary 5.5)} such that
 $(\al_j\,|\,\lam)>0$ for all $j\in{\bd J}\setminus{\bd K}$
 and $(\al_k\,|\,\lam)=0$ for all $k\in{\bd K}$.
 Choose a reduced expression
 $t_{\lam}\ft{|_{V_{\tiJ}}}=\bs\ft{(1)}\cdots\bs\ft{(n)}$ with
 $n\in{\mathbb N}$ and $\bs\ft{(1)},\dots,\bs\ft{(n)}\in S_{\tiJ}$,
 and set $\bs\ft{(p)}=\bs\ft{(\ol{p})}$ for each $p\in{\mathbb N}$,
 where $\ol{p}$ is the unique positive integer such that
 $1\leq\ol{p}\leq n$ and $\ol{p}\equiv p\mod n$.
 Then the infinite sequence $\bs=(\bs\ft{(p)})_{p\in{\mathbb N}}$
 is an element of ${\cl W}_{\tiJ}^\infty$ such that
 $\vPhi_{\tiJ}^\infty([\bs])=\De^{\tiK}_{\tiJ}\ft{(1,-)}$.
\end{prop}

\noindent
{\it Remark.}
 In \cite{jB}, the above proposition is proved in the case where
 ${\bd J}={\SI}$ and ${\bd K}=\emptyset$.

\begin{defn} 
 For each proper subset ${\bd K}$ of ${\bd J}$, we denote by
 $\ft{Z}^{\tiK}_{\tiJ}$ the unique element of $W_{\tiJ}^\infty$
 such that
 $\vPhi_{\tiJ}^\infty(\ft{Z}^{\tiK}_{\tiJ})=
 \De^{\tiK}_{\tiJ}\ft{(1,-)}$, and define a mapping
 $\chi_{\tiJ}\colon\bol{\cl P}_{\tiJ}\to W_{\tiJ}^\infty$
 by setting
\[ \chi_{\tiJ}(\ft{(\tiK,u,y)}):=uy.\ft{Z}^{\tiK}_{\tiJ} \]
 for each $(\scK,u,y)\in\bol{\cl P}_{\tiJ}$.
\end{defn}

\begin{lem} 
 For each ${\bd K}\subset{\bd J}$ and $y\in{\sW}_{\tiK}T_{\tiJ}$,
 we have
\begin{align}
 & y\De^{\tiK}_{\tiJ}\ft{(1,-)}\doteq\De^{\tiK}_{\tiJ}\ft{(1,-)}, \\
 & y\De^{\tiK}_{\tiJ}\ft{(1,-)}\setminus\Ome\subset
 \De^{\tiK}_{\tiJ}\ft{(1,-)}, \\
 & \{\vPhi_{\tiJ}(y)\setminus(-\Ome)\}\cap\re_{\tiK+}=
 \vPhi_{\tiJ}(y)\cap\re_{\tiK+},
\end{align}
 where $\Ome=y\De^{\tiK}_{\tiJ}\ft{(1,-)}\cap\re_{\tiJ-}$.
\end{lem}

\begin{proof}
 Since $\ol{y}\in{\sW}_{\tiK}$ we have
 $\ol{y}{\sDe}{}^{\tiK}_{\tiJ-}={\sDe}{}^{\tiK}_{\tiJ-}$
 by Lemma 4.1(3). Hence (7.1) follows from (iii) of Lemma 5.2(1).
 Moreover, by (ii) of Lemma 5.2(1), we have
\begin{equation}
 \ol{y\De^{\tiK}_{\tiJ}\ft{(1,-)}}\subset{\sDe}{}^{\tiK}_{\tiJ-}.
\end{equation}
 By the definition of $\Ome$, we have
\[ y\De^{\tiK}_{\tiJ}\ft{(1,-)}\setminus\Ome=
 y\De^{\tiK}_{\tiJ}\ft{(1,-)}\cap\re_{\tiJ+}. \]
 Thus (7.2) follows from (7.4).
 Moreover, by (7.4) we have
\[ -\Ome=(-y\De^{\tiK}_{\tiJ}\ft{(1,-)})\cap\re_{\tiJ+}
 \subset\De^{\tiK}_{\tiJ}\ft{(1,+)}, \]
 and hence $(-\Ome)\cap\re_{\tiK+}=\emptyset$.
 Thus we get
\[ \{\vPhi_{\tiJ}(y)\setminus(-\Ome)\}\cap\re_{\tiK+}=
 \vPhi_{\tiJ}(y)\cap\re_{\tiK+}. \qed \]
\renewcommand{\qed}{}
\end{proof}

\begin{thm} 
 Let ${\bd J}$ be an arbitrary non-empty subset of ${\SI}$.

 {\em(1)} For each $x\in W_{\tiJ}$ and ${\bd K}\subsetneq{\bd J}$,
 we have the following equality\/{\em:}
\begin{equation} 
 \vPhi_{\tiJ}^\infty(x.\ft{Z}^{\tiK}_{\tiJ})=
 \nab_{\tiJ}\ft{(\tiK,\ol{x}^{\tiK},z_x)}
\end{equation}
 with a unique element $z_x\in W_{\tiK}$ such that
\begin{equation} 
 \vPhi_{\tiJ}((\ol{x}^{\tiK})^{-1}x)\cap\De_{\tiK+}=
 \vPhi_{\tiK}(z_x).
\end{equation}

 {\em(2)} Both $\vPhi_{\tiJ}^\infty$ and $\chi_{\tiJ}$ are
 bijective and the following diagram is commutative\/{\em:}

\setlength{\unitlength}{0.5mm}
\begin{picture}(240,50)(-100,-5)
\put(10,30){\makebox(20,10){${\fr B}_{\tiJ}^\infty$}}
\put(-20,0){\makebox(20,10){$W_{\tiJ}^\infty$}}
\put(42,0){\makebox(20,10){$\bol{\cl P}_{\tiJ}$\;.}}
\put(45,5){\vector(-1,0){48}}
\put(-7,10){\vector(1,1){20}}
\put(47,10){\vector(-1,1){20}}
\put(17,0){$\chi_{\tiJ}$}
\put(-7,20){$\vPhi_{\tiJ}^\infty$}
\put(40,20){$\nab_{\tiJ}$}
%
\end{picture}

 {\em(3)} We have the following orbit decomposition\/{\em:}
\[ W_{\tiJ}^\infty=
 \coprod_{{\bd K}\subsetneq{\bd J}}W_{\tiJ}.\ft{Z}^{\tiK}_{\tiJ}. \]
\end{thm}

\begin{proof}
 (1) Put $y=(\ol{x}^{\tiK})^{-1}x$. Then $y\in{\sW}_{\tiK}T_{\tiJ}$.
 By Proposition 2.15, we have
\begin{equation} 
 \vPhi_{\tiJ}^\infty(y.\ft{Z}^{\tiK}_{\tiJ})=
 \{\vPhi_{\tiJ}(y)\setminus(-\Ome)\}\amalg
 \{y\De^{\tiK}_{\tiJ}\ft{(1,-)}\setminus\Ome\},
\end{equation}
 where $\Ome=y\De^{\tiK}_{\tiJ}\ft{(1,-)}\cap\re_{\tiJ-}$.
 Since $\sharp\Ome<\infty$ we have
 $\vPhi_{\tiJ}^\infty(y.\ft{Z}^{\tiK}_{\tiJ})
 \doteq\De^{\tiK}_{\tiJ}\ft{(1,-)}$
 by (7.1). Thus, by Lemma 6.1(2) we get
\begin{align} 
 & \De^{\tiK}_{\tiJ}\ft{(1,-)}\subset
 \vPhi_{\tiJ}^\infty(y.\ft{Z}^{\tiK}_{\tiJ}), \\
 & \De^{\tiK}_{\tiJ}\ft{(1,+)}\cap
 \vPhi_{\tiJ}^\infty(y.\ft{Z}^{\tiK}_{\tiJ})=\emptyset.
\end{align}
 By (7.2), (7.3), and (7.7), we have
\begin{equation} 
 \vPhi_{\tiJ}^\infty(y.\ft{Z}^{\tiK}_{\tiJ})\cap\re_{\tiK+}=
 \vPhi_{\tiJ}(y)\cap\re_{\tiK+}=\vPhi_{\tiK}(z_x),
\end{equation}
 where the second equality follows from (7.6).
 By (7.8)--(7.10) with (5.7), we have
\begin{align*}
 \vPhi_{\tiJ}^\infty(y.\ft{Z}^{\tiK}_{\tiJ})&=
 \{\vPhi_{\tiJ}^\infty(y.\ft{Z}^{\tiK}_{\tiJ})\cap
 \De^{\tiJ}\ft{(1,-)}\}\amalg
 \{\vPhi_{\tiJ}^\infty(y.\ft{Z}^{\tiK}_{\tiJ})\cap\re_{\tiK+}\} \\
 &= \De^{\tiK}_{\tiJ}\ft{(1,-)}\amalg\vPhi_{\tiK}(z_x)
 =\nab_{\tiJ}\ft{(\tiK,1,z_x)}.
\end{align*}
 Hence, by Proposition 2.15 and (6.2), we get
\[ \vPhi_{\tiJ}^\infty(x.\ft{Z}^{\tiK}_{\tiJ})=
 \vPhi_{\tiJ}^\infty(\ol{x}^{\tiK}.y.\ft{Z}^{\tiK}_{\tiJ})=
 \vPhi(\ol{x}^{\tiK})\amalg
 \ol{x}^{\tiK}\vPhi_{\tiJ}^\infty(y.\ft{Z}^{\tiK}_{\tiJ})=
 \nab_{\tiJ}\ft{(\tiK,\ol{x}^{\tiK},z_x)}. \]

 (2) By (1), we have
\begin{equation} 
 \vPhi_{\tiJ}^\infty(uy.\ft{Z}^{\tiK}_{\tiJ})=
 \nab_{\tiJ}\ft{(\tiK,u,y)}
\end{equation}
 for each $(\scK,u,y)\in\bol{\cl P}_{\tiJ}$. Hence
 $\vPhi_{\tiJ}^\infty\circ\chi_{\tiJ}=\nab_{\tiJ}$, which implies
 the surjectivity of $\vPhi_{\tiJ}^\infty$ since
 $\nab_{\tiJ}$ is bijective (see Theorem 6.7).
 Moreover, since $\vPhi_{\tiJ}^\infty$ is injective,
 $\vPhi_{\tiJ}^\infty$ is bijective, so is $\chi_{\tiJ}$.

 (3) Since $\chi_{\tiJ}$ is surjective, we have
 $W_{\tiJ}^\infty=
 \cup_{{\bd K}\subsetneq{\bd J}}W_{\tiJ}.\ft{Z}^{\tiK}_{\tiJ}$.
 Hence, it suffices to show that this union is disjoint.
 By (7.5), (7.11), and the injectivity of $\vPhi_{\tiJ}^\infty$,
 we have the following equality:
\begin{equation} 
 x.\ft{Z}^{\tiK}_{\tiJ}=\ol{x}^{\tiK}z_x.\ft{Z}^{\tiK}_{\tiJ}
\end{equation}
 for each $x\in W_{\tiJ}$. Suppose that
 $x.\ft{Z}^{\tiK}_{\tiJ}=y.\ft{Z}^{\tiL}_{\tiJ}$
 for some ${\bd L}\subsetneq{\bd J}$ and $y\in W_{\tiJ}$.
 By (7.12), we have
 $\ol{x}^{\tiK}z_x.\ft{Z}^{\tiK}_{\tiJ}=
 \ol{y}^{\tiL}z_{y}.\ft{Z}^{\tiL}_{\tiJ}$
 with a unique $z_{y}\in W_{\tiL}$.
 Thus we get ${\bd K}={\bd L}$ since $\chi_{\tiJ}$ is injective.
\end{proof}

\noindent
{\it Remark.}
 The existence and uniqueness of
 the element $z_x\in W_{\tiK}$ satisfying (7.6)
 are guaranteed by Lemma 2.5(2) and Corollary 5.5(1).

\begin{lem} 
 Let $B$ be a biconvex set in $\De_{\tiJ+}$.
 Then either $B\subset\re_{\tiJ+}$ or $\im_+\subset B$ holds.
\end{lem}

\begin{proof}
 We claim that if $B\cap\im_+\neq\emptyset$ then $\im_+\subset B$.
 Indeed, if $m\de\in B$ for some $m\in{\mathbb N}$, then we have
 $\de\in B$ by the convexity of $\De_{\tiJ+}\setminus B$,
 and hence $m\de\in B$ for all $m\in{\mathbb N}$
 by the convexity of $B$, i.e., $\im_+\subset B$.
 Thus either $B\subset\re_{\tiJ+}$ or $\im_+\subset B$ holds.
\end{proof}

\begin{cor} 
 Let $B$ be a subset of $\De_{\tiJ+}$.
 Then $B$ is a biconvex set in $\De_{\tiJ+}$ if and only if
 one of the following {\em (a)--(d)} holds\/{\em:}
\[ \text{\em(a)}\; B=\vPhi_{\tiJ}(z); \quad
 \text{\em(b)}\; B=\De_{\tiJ+}\setminus\vPhi_{\tiJ}(z);\quad
 \text{\em(c)}\; B=\vPhi_{\tiJ}^\infty(\ft{Z});\quad
 \text{\em(d)}\; B=\De_{\tiJ+}\setminus\vPhi_{\tiJ}^\infty(\ft{Z}), \]
 where $z$ is an element of $W_{\tiJ}$ and
 $\ft{Z}$ is an element of $W_{\tiJ}^\infty$.
\end{cor}

\begin{proof}
 The ``if'' part is obvious.
 Let us prove the ``only if'' part.
 By Lemma 7.5, we have either
 $B\subset\re_{\tiJ+}$ or $\im_+\subset B$.
 If $B\subset\re_{\tiJ+}$ and $\sharp B<\infty$, then
 $B=\vPhi_{\tiJ}(z)$ with $z\in W_{\tiJ}$ by Corollary 5.5(1).
 If $B\subset\re_{\tiJ+}$ and $\sharp B=\infty$, then 
 $B=\vPhi_{\tiJ}^\infty(\ft{Z})$ with $\ft{Z}\in W_{\tiJ}^\infty$
 by Theorem 7.4(2).
 If $B\subset\pim$, then $\De_{\tiJ+}\setminus B$
 is a real biconvex set in $\De_{\tiJ+}$.
 Hence we have either
 $\De_{\tiJ+}\setminus B=\vPhi_{\tiJ}(z)$ or
 $\De_{\tiJ+}\setminus B=\vPhi_{\tiJ}^\infty(\ft{Z})$,
 i.e.,
 $B=\De_{\tiJ+}\setminus\vPhi_{\tiJ}(z)$ or
 $B=\De_{\tiJ+}\setminus\vPhi_{\tiJ}^\infty(\ft{Z})$,
 where $z\in W_{\tiJ}$ and $\ft{Z}\in W_{\tiJ}^\infty$.
\end{proof}

\noindent
{\it Remark.}
 By the corollary, we see that a subset $B\subset\De_{\tiJ+}$
 is a biconvex set in $\De_{\tiJ+}$ if and only if $B$ satisfies
 the conditions $\mbox{C(i)}'$ and $\mbox{C(ii)}'$
 with replacing $\De_+$ by $\De_{\tiJ+}$
 (see the remarks below Theorem 2.6).
